\documentclass[smallcondensed]{svjour3}     
\usepackage{listings}
\usepackage{color}
\title{Syntax Highlighting in LaTeX with the listings Package}
\author{writeLaTeX}
\definecolor{mygreen}{rgb}{.9,1,0.7}
\definecolor{mygray}{rgb}{0.5,0.5,0.5}
\definecolor{mymauve}{rgb}{0.58,0,0.82}
\smartqed  
\usepackage{amsmath, amscd, amsfonts, amssymb, graphicx}
\usepackage[a4paper]{geometry}
\usepackage{siunitx}
\usepackage{mathrsfs}
\newtheorem{thm}{Theorem}[section]

\newtheorem{defn}{Definition}[section]

\usepackage{hyperref}
\usepackage{cite}
\usepackage{soul,color}
\usepackage[english]{babel}
\usepackage[utf8x]{inputenc}
\usepackage{fullpage}
\begin{document}
\title{Bifurcation analysis of a fractional-order Pinsky-Rinzel model
}
\titlerunning{Fractional-order CA3 hippocampal pyramidal neurons}       
\author{}

\author{Leila Eftekhari\and Soleiman Hosseinpour \and Moein Khalighi \and Salvador Jim\'{e}nez 
}
\institute{Leila Eftekhari \at 
Department of Mathematics,  Tarbiat Modares University, Tehran, Iran\\
\email{leila.eftekhari32@gmail.com} 
\and
Soleiman Hosseinpour \at
Department of Applied Mathematics, Shahrood University of Technology, Shahrood, Iran\\
\email{soleiman.hosseinpour@gmail.com} 
 \and
 Salvador Jim\'{e}nez \at
Departamento de Matem\'{a}tica Aplicada a las TIC,
Universidad Polit\'{e}cnica de Madrid, Madrid, Spain\\
\email{s.jimenez@upm.es}
\and
Moein Khalighi \at
Department of Computing, University of Turku, Finland\\
 \email{moein.khalighi@utu.fi}
 }
\date{Received: date / Accepted: date}
\maketitle
\begin{abstract}
The present work describes a new fractional-order system of a two-compartment CA3 hippocampal pyramidal cell, so-called the Pinsky-Rinzel model, with Caputo fractional derivative. We investigate the transient of the solutions and the general behavior of the system based on the bifurcation diagrams, in which fractional derivative orders and currents injection are taken as bifurcation parameters. Chaotic regions are obtained by tunning the values of the fractional derivative orders and injection currents. The proposed numerical approach helps to study the stability of the system under certain conditions.

\keywords{CA3 hippocampal pyramidal neurons, Caputo fractional derivative, Bifurcation analysis.}
\end{abstract}

\section{Introduction}
The CA3 region in the brain system has received more attention in recent years due to its particular role in memory encoding \cite{kandel2000principles}. This area is critical for initiating hippocampal interictal-like activity, which depends on the presence of abundant recurrent excitatory connections between CA3 pyramidal cells \cite{miles1983single }. CA3 pyramidal cells intrinsically display various firing patterns at the cellular level, ranging from single action potentials to complex bursts \cite{bib3}.  Such bursting is important for place-cell activity signal propagation and the induction of synaptic plasticity \cite{bib4}. On the other hand, mathematical models aim to understand the behavior of CA3 pyramidal cells \cite{bib5}. In 1994, Pinsky and Rinzel developed a two-compartment model for CA3 hippocampal pyramidal neurons in a guinea pig based on a complex 19 compartments Traub model  \cite{bib1}. Their proposed model explains how interactions between the somatic and dendritic compartments occur, thus allowing a computational implementation and attracting many scholars' attention \cite{bib4}.

The Pinsky-Rinzel model is a non-smooth, and classic methods' mathematical analysis of its dynamical behavior is inefficient. In 2001, Hahn and Drand carried out a mathematical analysis of the dynamical properties of the model based on Strogatz's works\cite{bib3}. They probed the changeover between resting, bursting, and spiking states affected by increasing the amount of extracellular potassium concentration. Due to the complexity and non-smooth nature of the system, analyzing the dynamic scenery of the system was quite difficult. For these reasons, L. A. Atherton and \textit{et al.} modified the Pinsky-Rinzel model \cite{bib3}. 
These models investigate both qualitative acontinuous functionsnd quantitative analysis of the  CA3 pyramidal cells, ranging from detailed, multi-compartmental models down to single compartmental models. 
This is important since it allows performing bifurcation analysis and stability investigations and paves the way to implementing more numerical methods. 

A bifurcation diagram is a good representation to identify the disappearance, appearance, or stability changes like the limit sets (equilibrium points, periodic or almost periodic solutions) as the bifurcation parameter varies. This concept has been significantly studied in various fields during the last decades due to its importance for understanding the behavior of dynamic systems \cite{PhysRevResearch.2.023281, article1, article2}. In the field of Neuroscience,  authors in \cite{bib3} reported a bifurcation analysis of a two-compartment hippocampal pyramidal cell model, which is the only bifurcation analysis of the CA3 pyramidal cell model until now. They used a smooth system to investigate the qualitative behavior of the original model. However, in the current work, we generalized a CA3 pyramidal cell model, known as the Pinsky-Rinzel model, with the concept of fractional calculus to study the behavior of the action optional and its impact on the gating variables in the proposed model. Fractional calculus is a generalization of ordinary derivatives and integrals to non-integer order. It has been shown that fractional calculus operators hold non-local properties and can represent memory effects. In dynamical systems theory, these effects are described by explicit non-local terms of past system states \cite{nagy2014efficient}.  

Many new fractional-order systems have been suggested in recent years that they have rich and complex dynamics, such as periodic solutions, chaos windows, and bifurcation. For instance, the Chua system with a derivative order 2.7 makes chaos motion \cite{hartley1995chaos}. Chaotic dynamics of a damped van der Pol equation with fractional order are reported in \cite{chen2008chaotic}. Authors in  \cite{gao2005chaos} investigated the nonlinear dynamics and chaos in a fractional-order financial system. A periodically forced complex Duffing’s oscillator was proposed, and chaos for the system was introduced in detail \cite{chen2008nonlinear}. Bifurcation, chaos control, and synchronization were extended for a fractional-order Lorenz system with complex variables \cite{gao2005chaos}. 
In fractional-order systems, a bifurcation implies a qualitative or topological change in dynamics with a variation of a system parameter or derivative order, and bifurcation analysis becomes harder due to the non-local property of the operator of fractional calculus. Concerning the importance of bifurcation analysis in fractional-order systems, the number of studies on fractional-order models is still limited in Neuroscience.

In this paper, motivated by the above discussion, we survey the solutions (membrane potentials and currents) of the fractional-order Pinsky-Rinzel model. We study the bifurcations of the model as the fractional derivative order changes, and we identify regular and chaotic regimes. Next, we consider the somatic and dendritic injections, ISapp, and IDapp , as bifurcation parameters and explore the chaotic behavior of the system in each case. Finally, we investigate the system's stability condition and try to find some conditions that make the equilibria asymptotically stable.

\begin{figure}
\centering
\includegraphics[width=\textwidth]{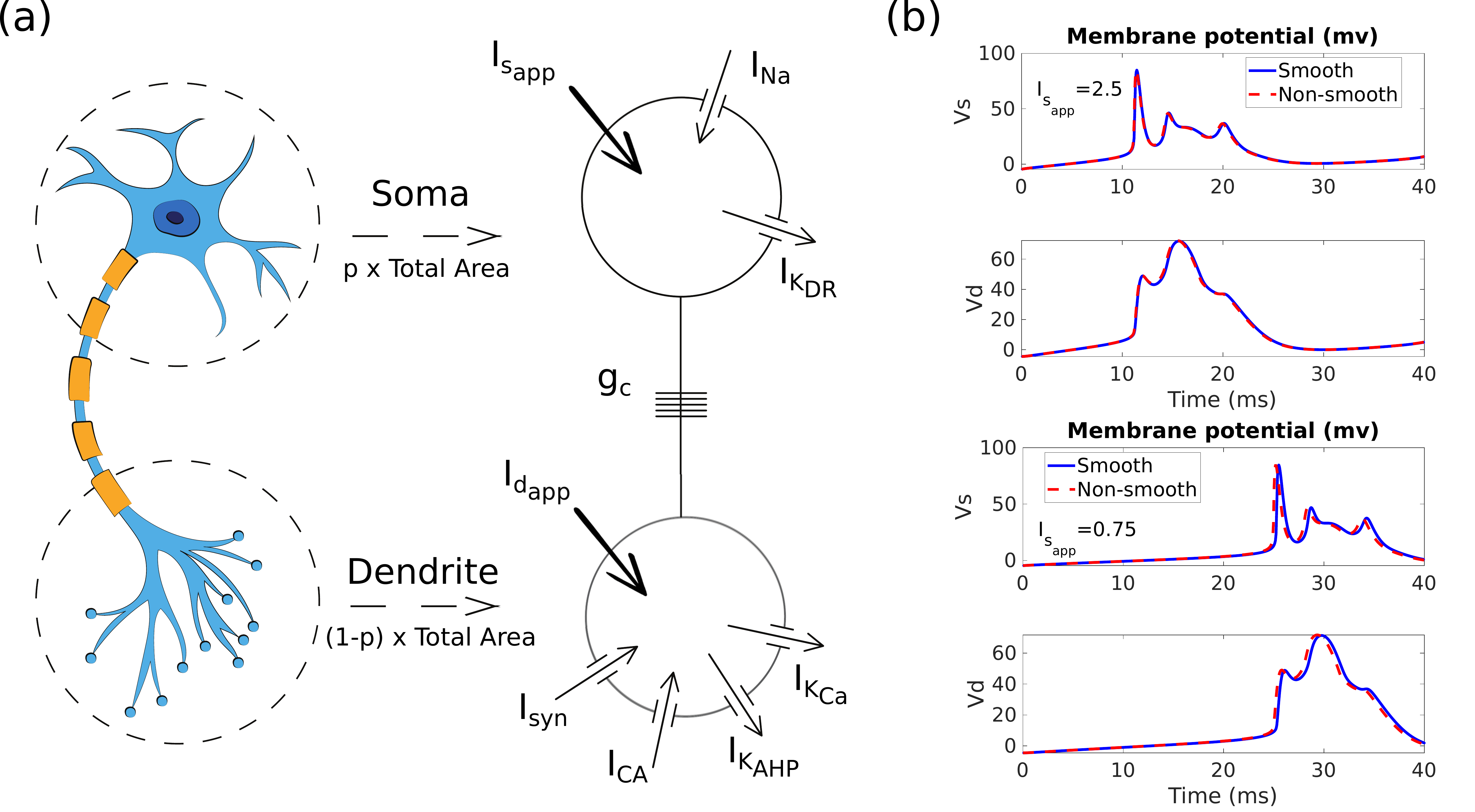}
\caption{\textbf{Schematic of the model and its dynamical behavior}. (\textbf{a}) A two-compartment model \eqref{ca3} displaying currents from and to soma and dendrite. The somatic compartment is combined with fast Sodium $\si{I_{Na}}$  and delayed rectifier Potassium $\si{I_{K_{DR}}}$ and leak current. The dendritic compartment has a persistent Calcium $ \si{I_{Ca}}$, Calcium activated Potassium $ \si{I_{K_{Ca}}}$, hyperpolarisation Potassium current  $\si{I_{K_{AHP}}}$, synaptic current \(I_{syn}\), and leak current. Leak currents have been omitted for clarity. Electronic coupling between  the two compartments is modeled using two parameters; \si{g_{c}} is  the strength of coupling, and  $p$ is the percentage of the total area in the somatic-like compartment. $\si{I_{S_{app}}}$ and $\si{I_{D_{app}}}$ are currents injected through an electrode to couple the currents between the two compartments. (\textbf{b}) Non-smooth  curve and smooth approximation of $\si{V_s}$  and  $\si{V_d}$ from  the model \eqref{ca3} for two particular values of $\si{I_{S_{app}}}=2.5$ and 0.75 when $\si{I_{D_{app}}}=0$.
When the membrane potential is close to the resting potential of the cell, these channels stay in a closed state, but when the membrane potential increases to a certain threshold value, these channels start to open quickly. At this moment, the action potential is fired. When the action potential dies, the membrane enters a hyperpolarized state; during this stage, no action potential can fire. This period is considered  the refractory period.
 }\label{fig1}
\end{figure}
\newpage
\section{Model} 
The following sections describe the main features of the CA3 hippocampal pyramidal neurons modeled by an ordinary differential equation and show how to incorporate non-integer order derivatives into this model.
\subsection{Pinsky-Rinzel model}
The Pinsky-Rinzel model is based on two compartments illustrated in Fig.~\ref{fig1}a. The CA3 mathematical model proposed in \cite{bib3, bib4, bib5} is described by an integer-order system of differential equations 
\begin{equation}
\begin{aligned}\label{ca3}
&\si {C_m}\si {V_s^{'}}=-\si {I_{Leak}}-\si {I_{Na}}-\si {I_{K_{DR}}}-\dfrac{\si {I_{DS}}}{\si p}+\dfrac{\si {I_{S_{app}}}}{\si p},\\
&\si {C_m}\si {V_d^{'}}=-\si{ I_{Leak}}-\si {I_{Ca}}-\si {I_{K_{Ca}}}-\si {I_{K_{AHP}}}+\dfrac{\si{ I_{SD}}}{(1-\si p)}+\dfrac{\si {I_{D_{app}}}}{(1-\si p)}.
\end{aligned}
\end{equation}
The currents of the model are functions of potentials as follows
\begin{align*}
\si{I_{Na}}&=\si{g_{Na}}{m_\infty^2}(\si{V_s})h(\si{V_s}-\si{V_{Na}}),\\
\si{I_{K_{DR}}}&=\si{g_{K_{DR}}} n(\si{V_s}-\si{V_{K}}),\\
\si {I_{Ca}}&=\si{g_{Ca}s}^2(\si{V_d}-\si{V_{N}}),\\
\si {I_{K_{Ca}}}&=\si{g_{k_{Ca}}} C~ \chi(\si Ca)(\si{V_d}-\si{V_{Ca}}),\\
\si{I_{K_{AHP}}}&=\si{g_{K_{AHP}}} q~(\si{V_d}-\si{V_{K}}),\\
\si {I_{SD}}&=-\si{I_{DS}}=\si{g_{c}}(\si{V_d}-\si{V_{s}}),\\
\si {I_{Leak}}&=\si{g_{L}}(\si{V}-\si{V_{L}}).
\end{align*}
The activation and  inactivation variables should satisfy these equations
\begin{align}
\omega^{'}(\si V)&=\dfrac{\omega_\infty(\si V)-\omega}{\tau_\omega(\si V)},\\
\omega_\infty(\si V)&=\dfrac{\alpha_\omega(\si V)}{\alpha_\omega(\si V)+\beta_\omega(\si V)},\\
\tau_\omega(\si V)&=\dfrac{1}{\alpha_\omega(\si V)+\beta_\omega(\si V)},
\end{align}
where, independly, we consider $\omega=h, n, s, m, C$ and $ q$. The rate functions are defined as follows
\begin{align*}
\alpha_m(\si {V_s})&=\dfrac{0.32(-46.9-\si{V_s})}{\exp(\frac{-46.9-\si{V_s}}{4})-1},\\
\beta_m(\si{V_s})&=\dfrac{0.28(\si{V_s}+19.9)}{\exp(\frac{\si{V_s}+19.9}{5})-1},\\
\alpha_n(\si{V_s})&=\dfrac{0.016(-24.9-\si{V_s})}{\exp(\frac{-24.9-\si{V_s}}{5})-1},\\
\beta_n(\si{V_s})&=0.25\exp(-1-0.025\si{V_s}),\\
\alpha_h(\si{V_s})&=0.128\exp\Big(\dfrac{-43-\si{V_s}}{18}\Big),\\
\beta_h(\si{V_s})&=\dfrac{4}{1+\exp(\frac{(-20-\si{V_s})}{5})},\\
\alpha_s(\si{V_d})&=\dfrac{1.6}{1+\exp(-0.072(\si{V_d}-5))},\\
\beta_s(\si{V_d})&=\dfrac{0.02(\si{V_d}+8.9)}{\exp(\frac{(\si{V_d}+8.9)}{5})-1},\\
\alpha_C(\si{V_d})&=\dfrac{(1-H(\si{V_d}+10))\exp(\frac{(\si{V_d}+50)}{11}-\frac{(\si{V_d}+53.5)}{27})}{18.975}+H(\si{V_d}+10)(2\exp(\frac{(-53.5-\si{V_d})}{27})),\\
\beta_C(\si{V_d})&=(1-H(\si{V_d}+10))(2\exp(\frac{(-53.5-\si{V_d})}{27})-\alpha_c(\si{V_d})),\\
\alpha_q(\si{Ca})&=\min(0.00002\si{Ca}, 0.01),\\
\beta_q(\si{Ca})&=0.001,\\
{\chi}(\si{Ca})&=\min\Big(\dfrac{\si{Ca}}{250}, 1\Big),
\end{align*}
where $ w=h, s, n, m $ are defined simply by continuous rate functions, while  the $C$,  $q$, and $\chi$ are formulated as discontinuous rate functions where $H(.)$ is the Heaviside step function. However, this is a non-smooth approximation, and for bifurcation analysis, these discontinuous functions can be converted to a smooth form \cite{bib3}:
\begin{align*}
C_\infty(\si{V_d})=&\left(\dfrac{1}{1 + \exp(\dfrac{-10.1 - V_d}{0.1016})}\right) ^ {0.00925},\\
\tau_C(\si{V_d})=&3.627  \exp(0.03704 V_d),\\
q_\infty(\si{Ca})=&0.7894 \exp(0.0002726  \si{Ca})-0.7292 ~\exp(-0.01672 \si{ Ca}),\\
\tau_q(\si{Ca})=&657.9 \exp(-0.02023 \si{ Ca})+301.8 \exp(-0.002381 \si{ Ca}),\\
\chi(\si{Ca})=&1.073  \sin(0.003453 \si{ Ca} + 0.08095)+0.08408 \sin(0.01634 \si{Ca} - 2.34)\\
+&0.01811\sin(0.0348 ~ \si{Ca} - 0.9918),
\end{align*}
also we have 
\begin{align*}
\dfrac {\si{dCa}}{\si{dt}}= -0.13\si I_{\si{Ca}}-0.075\si{Ca}.
\end{align*}
The first minus sign is based on the convention that inward currents are negative. The non-smooth and the smooth approximation demonstrate very similar dynamic behaviors (Fig. \ref{fig1}b). In Appendix, the notations and maximal conductance parameters are given in  Table \ref{t2}. 
\subsection{Fractional-order formulation}
In this section, a time-dependent Pinsky-Rinzel model is introduced by using fractional calculus. There are several definitions for fractional derivatives and integrals for different purposes~\cite{de2014review}. The Caputo fractional derivatives are more applicable due to providing initial conditions with clear physical interpretation. It motivates us to propose the CA3 model in the sense of the Caputo fractional time derivative, which holds non-local properties. 

The Caputo fractional time derivative of order $ \alpha $ is defined by  \cite{bib8}
\begin{align}
{}^{C}\!D^\alpha f(t)=
\begin{cases}
I^{n-\alpha}D^n f(t),&~~n-1<\alpha <n,\\
D^n f(t),&~~\alpha=n.
\end{cases}
\end{align}
where $ n\in \mathbb{N}$, $D^n$ is the integer derivative of order $n$, and \(I^{n-\alpha}\) is the factional integral of order $(n-\alpha)$ over function $ f(t) $ that is defined by
\begin{align}
I^{n-\alpha} f(t)=\dfrac{1}{\Gamma(n-\alpha)}\int_0^t f(\tau)(t-\tau)^{n-1-\alpha}d\tau.
\end{align}
where $\Gamma(.)$ is Euler's Gamma function and $f(t)\in C^n$.
Replacing the integer order derivatives of \eqref{fig1} by non-integer order derivatives causes dimensional inconsistencies in the equation. Hence, we should insert some parameters in the fractional order model to make the same unit of measure for both sides of the equations. Westerlund and Ekstam \cite{bib26} have shown that Jacques Curie's empirical law for the current through capacitors and dielectrics causes a capacitive current-voltage relationship for a non-ideal capacitor
 
\[I_{C}^\alpha=C_m(\alpha)\,{}^C\!D^\alpha\]
where $\alpha \in (0, 1)$, $C_{m}(\alpha)$ indicates the fractional-order capacitance with units $\si{(A/V)s^{\alpha} }$ and $ \,{}^C\!D^\alpha$ is the Caputo fractional-order differential operator. 
Here, we generalize the Pinsky-Rinzel model~\eqref{ca3} by considering fractional-order derivatives of the somatic  ($\si{V_s}$) and dendritic ($\si{V_d}$) membrance potentials
\begin{equation} 
\begin{aligned}
&C_m(\alpha)\,{}^C\!D^\alpha  V_s=-\si {I_{Leak}}-\si {I_{Na}}-\si {I_{K_{DR}}}-\dfrac{\si {I_{DS}}}{\si p}+\dfrac{\si {I_{S_{app}}}}{\si p},\\
&C_m(\alpha)\,{}^C\!D^\alpha  V_d=-\si{ I_{Leak}}-\si {I_{Ca}}-\si {I_{K_{Ca}}}-\si {I_{K_{AHP}}}+\dfrac{\si{ I_{SD}}}{(1-\si p)}+\dfrac{\si {I_{D_{app}}}}{(1-\si p)}.
\end{aligned}\label{bif1}
\end{equation}
where $C_m(\alpha)=\dfrac {\tau^\alpha}{R_m}$, $R_{m}$ is the membrane resistance, and $ \tau$ is the time constant \cite{bib27}.
Also we rewrite 
\begin{align}
\,{}^C\!D^\alpha\omega(\si V)&=\dfrac{\omega_\infty(\si V)-\omega}{\tau_\omega(\si V)},
\end{align}
where $ \omega={h, s, n, m }$, $C$, $q$ and $ {\chi}$  are defined simply  by continuous rate functions in the same way. 
\section{Results}
\begin{figure}
\centering
\includegraphics[width=\textwidth]{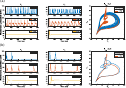}
\caption{ Time series plots of $\si{V_s}$ and  $\si{V_d}$ from model \eqref{ca3} for two particular values of $\si{I_{S_{app}}}=2.5$ and 0.75 when $\si{I_{D_{app}}}=0$ and time = 500.}\label{fig2}
\end{figure}
This section discusses some dynamical properties, transient, and possible chaotic behavior, that the fractional-order Pinsky-Rinzel reduced model modifies.  Our numerical method is based on the predictor-corrector Adams-Bashforth-Moulton method for fractional differential equations \cite{dithelm, bib28}.

The numerical solutions of the model is provided for two types of interest somatic current injection $\si{I_{S_{app}}}=2.5$ and 0.75 $(\si{\mu A/cm^2})$ (Fig.\ref{fig25cur}, \ref{fig075cur})
Results display somatic membrane potential in two distinctive firing modes. The higher currents produce regular spiking (Fig.~\ref{fig2}a), whilst low somatic current injections cause low-frequency bursts in the system (Fig.~\ref{fig2}b). 
In the first mode, the dynamic of variable $q$ governing activation of $\si{ I_{K_{AHP}}}$ and calcium concentration stays at a high level to prevent a dendritic spike (Fig. \ref{fig25var}). In contrast, in the second mode, a burst is operated when the slow, hyperpolarizing current activation variable decreases below a threshold, causing a calcium spike in the dendrites (Fig. \ref{fig075var}). 

 The Action potential is produced in neurons , and the movement of ions does propagation through a special kind of voltage-gated ion channels in the membrane \cite{bib30}. when the membrane potential is close to the resting potential of the cell these channels stay in a closed state but when the membrane potential increases to a certain threshold value, these channels start to open quickly. At this moment, the action potential is fired. When the action potential dies, the membrane enters a hyperpolarized state; ; no action potential can fire during this stage. This period is considered as the refractory period.
\subsection{Studying transient}

Plotting the variables versus time can only show whether or not the solutions have reached the asymptotic regime. To comprehend of the proposed model and for future development in bifurcation analysis, we study transients of the system.
\begin{figure}
\centering
\includegraphics[width=15cm]{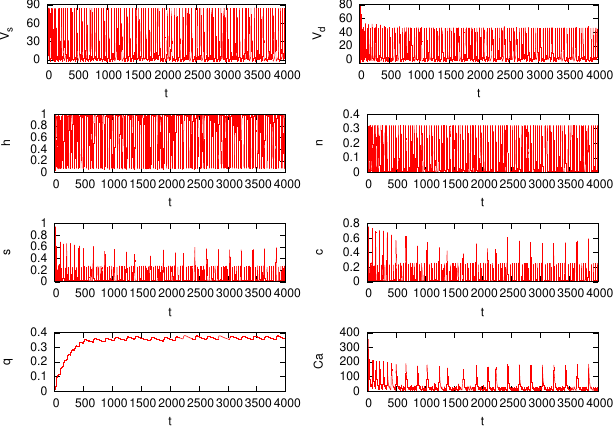}
\caption{ Time series plots of  $\si{V_s}$, $\si{V_d}$, Ca, h, q, n, s and c ~ from model \eqref{bif1} when  $\si{I_{S_{app}}}=2.5$, $\alpha=0.95$ and time=4000.}\label{fig3}
\end{figure}
Here, the transient is a regime in which the solution has a qualitative behavior quite different than what is observed afterward. Usually, this corresponds to an initial “short” time window after which the qualitative behavior of the system is more uniform and it is in opposition to a "steady state".

As can shown in (Fig. \ref{fig3}), with $\si{I_{S_{app}}}=2.5$,  variables $\si{V_s}$, $h$, $n$ have short transient. On the other hand, $\si{V_d}$, $s$, and Ca do not have short transient and have drastic amplitude changes for times up to $\sim 67$ (the distance between the time of the big peaks, which is an estimation). Finally, after $t=400$, the asymptotic regime seems to be obtained.  the asymptotic regime, follows the transient; the new qualitative regime of the solutions observed for longer times and that one may represent the behavior as t goes to infinity. 

\begin{figure}
\centering
\includegraphics[width=15cm]{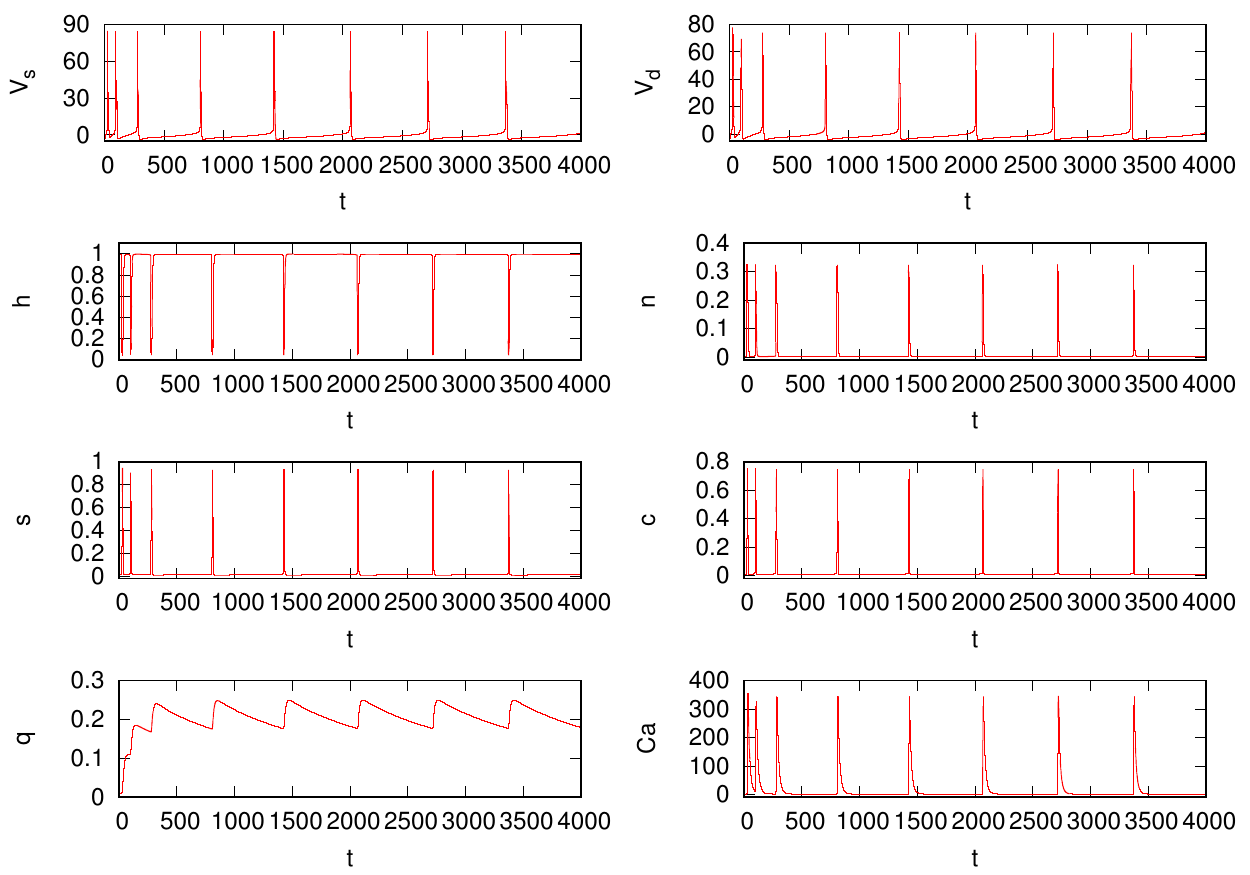}
\caption{ Time serious plots of  $\si{V_s}$, $\si{V_d}$, Ca, h, q, n, s and c ~ from model \eqref{bif1} when $\si{I_{S_{app}}}=0.75$, $\alpha=0.95$ and time=4000.}\label{fig4}
\end{figure}
Variables $c$ and $q$ do not present the initial burst, but they have a transient such that $t=400$, which means for these variables the behavior after 400 or even 500 is quite different from the one before, and this change is even much more evident for q (the average value grows and seem to reach a plateau after 400 or 500). Furthermore, the system seems to reach the asymptotic behavior after the transient, as can be seen if we extend the plot to higher times (an asymptotic behavior can be chaotic or regular). Variables $V_s$ and h do not show a change,  $V_d$ has an initial peak almost at double amplitude than the rest, and the six following peaks appear after shorter times than all the other peaks that appear later. Something similar occurs for variable n;  although the amplitude of the peaks is consistent in the graphic the times at which they occur are much shorter below 400 than afterward. For variable s, the initial peak is very high, the amplitude decreases and there is a change of behavior at about $t \approx 300$ and variable Ca also changes the behavior at about $t \approx 500$ ( Fig.\ref{fig3}).\\
\begin{figure}
\centering
\begin{tabular}{cc}
\includegraphics[width=6cm]{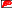}
&\includegraphics[width=6cm]{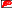}
\\
\includegraphics[width=6cm]{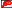}
&\includegraphics[width=6cm]{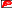}
\end{tabular}
 \caption{ $\si{V_d}$  versus $\si{V_s}$  model \eqref{bif1} when $\si{I_{S_{app}}}=2.5$ and  $\alpha=0.95$ }\label{figtime}
\end{figure}
\begin{figure}
\centering
\begin{tabular}{cc}
\includegraphics[width=6cm]{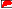}&
\includegraphics[width=6cm]{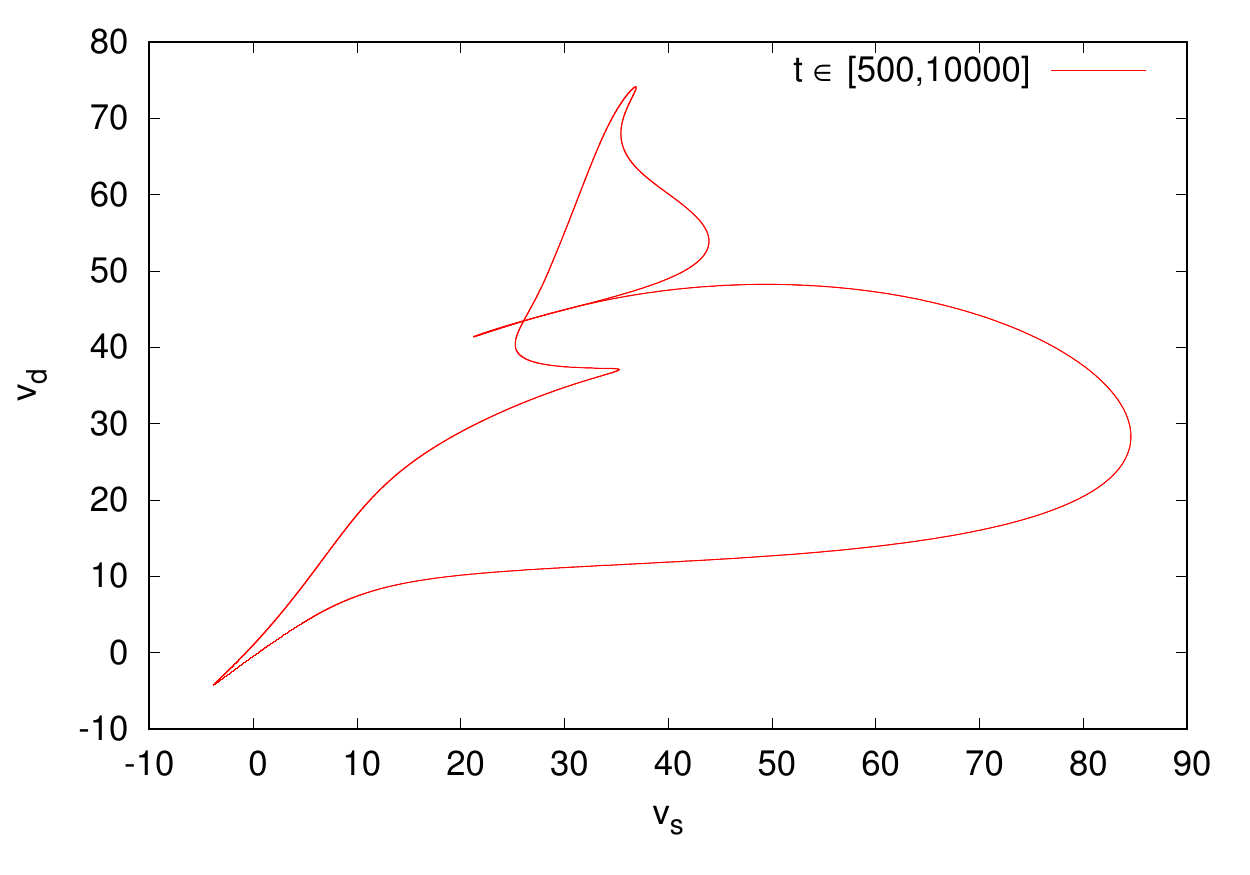}
\end{tabular}
\caption{ Attractor from  model \eqref{bif1} via $\si{V_s}$, $\si{V_d}$, $\alpha=0.95$, case  $\si{I_{S_{app}}}=2.5$ (left) and case $\si{I_{S_{app}}}=0.75$ (right) .}\label{GRAF}
\end{figure}
Furthermore, something similar appears in the case of $\si{I_{S_{app}}}=0.75$;  the transient does not show a difference in the amplitude but the periodical behavior. The transient for this case may be $[0,300]$. After that, the solutions seem to have a periodic behavior with period $\sim 67$ (see Fig. \ref{fig4}).
For more investigation, the final time is extended up to $t = 10000$.  The final time has extended drastically  because  a time-span $\Delta t=1000$ is not large enough to reveal the behaviour of this system since the data show that there is some intrinsic time-scale of the order of $t_s\approx 30$; $\Delta t\approx 33 ~ t_s$ is a little short. To make sure, the time initially is extended to $t = 2000$, then to $t = 3000$, and finally to $t = 10000$ (indeed, in the end, $t\approx 333~ t_s$). In order to simplify the study, some graphics of the first two variables are represented, $V_s$ and $V_d$ (see Fig. \ref{figtime}). \\
There is a transient (not very long; $t\approx400$) after which the solution is confined to values of $V_d$ in $[-5,50]$ (for $t\in[0,10000]$)  instead of $[-10,80]$ (for $t\in[1000, 2000]$), otherwise, the solution gives similar results. Since the number of lines is different for the case, $t\in[2000, 3000]$, to check that this does not go to a regular solution for longer times, the final time is extended to $t=10000$. Then, the behavior is the same; the system is chaotic (at least up to this time), and the final graphic for times $t\in[500,10000]$ shows what appears to be a strange attractor (Fig. \ref{GRAF}). If we do the same for the case, ${I_{S_{app}}}=0.75$, we see that the solution is  regular. If we remove the initial transient, we obtain a simple, closed curve (Fig. \ref{GRAF}). 
\subsection{Bifurcation analysis}
Bifurcations significantly aim to understand the dynamics of fractional-order systems \cite{bib24, bib25, bib12}, for instance, the Chen model \cite{bib24a} and the Hodgkin-Huxley model \cite{bib15}. The most remarkable point about bifurcation diagrams is that each branch looks like the whole diagram. The diagram contains many such copies of itself. This self-similarity can be seen to repeat itself at ever finer resolutions. Such behavior is characteristic of geometric entities called fractals \cite{devaney1990chaos}. \\

In this section, bifurcation analysis of system (\ref{bif1}) is carried out to investigate the dynamical behavior of the fractional-order Pinsky-Rinzel model \eqref{bif1} in terms of $\si{I_{S_{app}}}, \si{I_{D_{app}}}$, and order of fractional derivative $\alpha$.
 Qualitative changes in the system dynamics are called bifurcations, and the parameter values at which they occur are called bifurcation points.
We consider $\si{I_{S_{app}}}=2.5, \si{I_{D_{app}}}=0, \si{g_c}=2.1$ and  $\alpha\in I=(0,1)$. The fractional-order $\alpha$ is considered as the  bifurcation parameter. The bifurcation diagram demonstrates the chaotic behavior of the system. 
\begin{figure}[h!t] 
\centering 
\includegraphics[width=7cm,height=6cm]{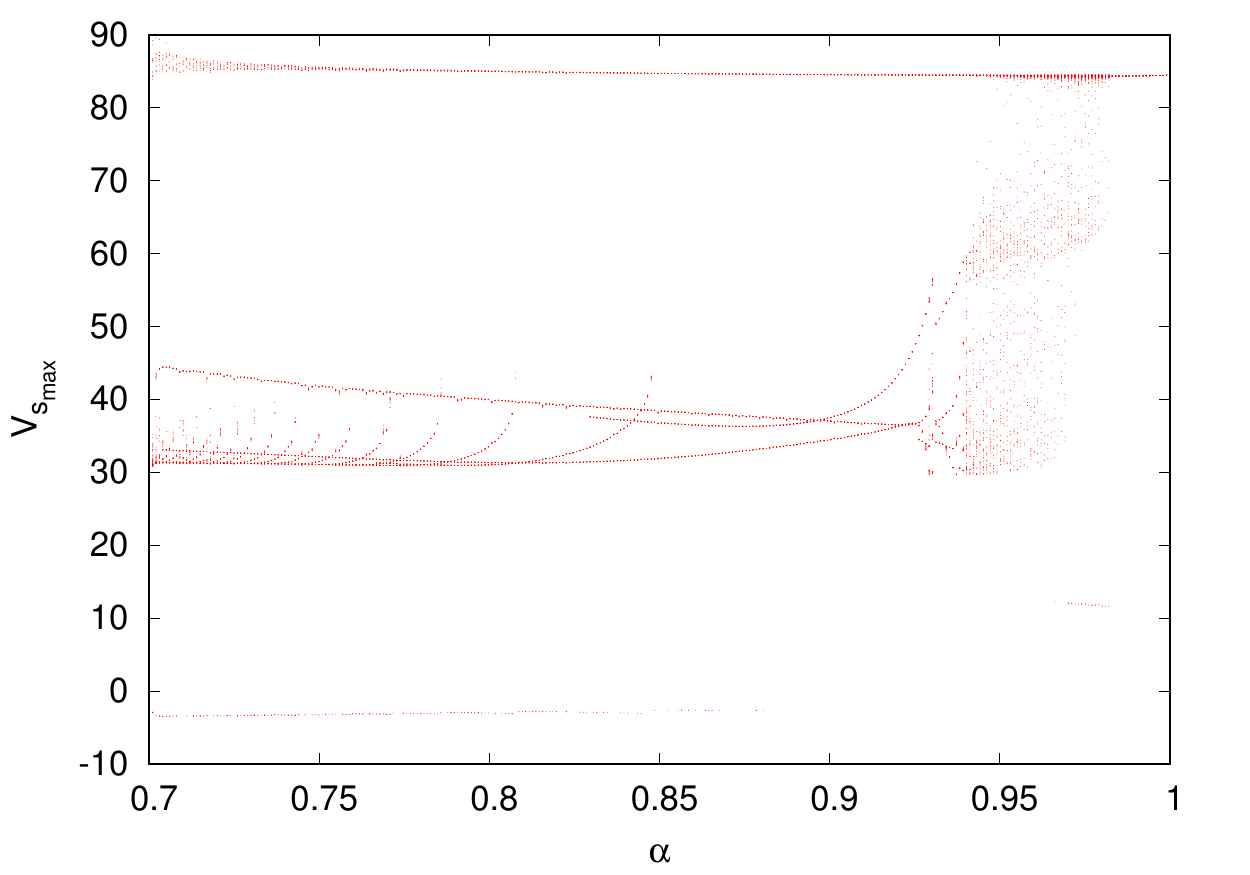}
\caption{ Bifurcation diagram for $\alpha\in[0.7,1]$ and $(\si{I_{S_{app}}}=2.5 $,  $ \si{I_{D_{app}}}=0,  $ and $ \si{g_c=2.1})$.}\label{fignew1} 
\end{figure}
According to Fig. \ref{fignew1}, the system is not chaotic for a range of $\alpha$, from 0.7 to \(\sim\)0.93. There is probably some transient behavior of the trajectories.
Subsequently, we turn $\si{I_{S_{app}}}$ as the bifurcation parameter when $ \si{I_{D_{app}}}=0, \si{g_c}=2.1$. Fig. \ref{fig11} shows that the order of derivatives \(\alpha\) alter the chaotic behaviors of the system. In fact, by increasing the $\si{I_{S_{app}}}$, chaos stops sooner in the systems with fewer order values.  
\begin{figure}[!tbp] 
\centering
\begin{tabular}{cc}
\includegraphics[width=7cm,height=6cm]{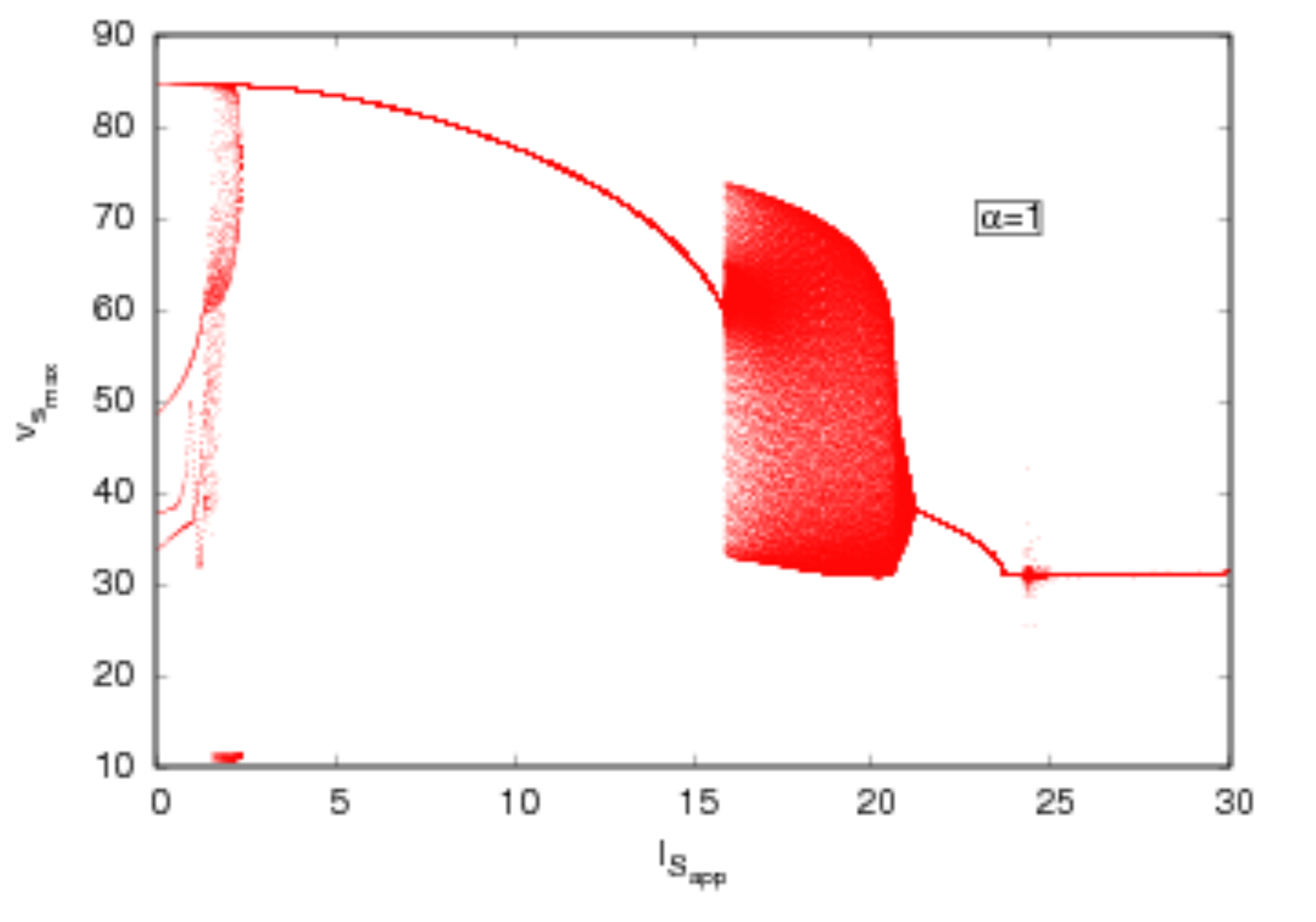}&
\includegraphics[width=7cm,height=6cm]{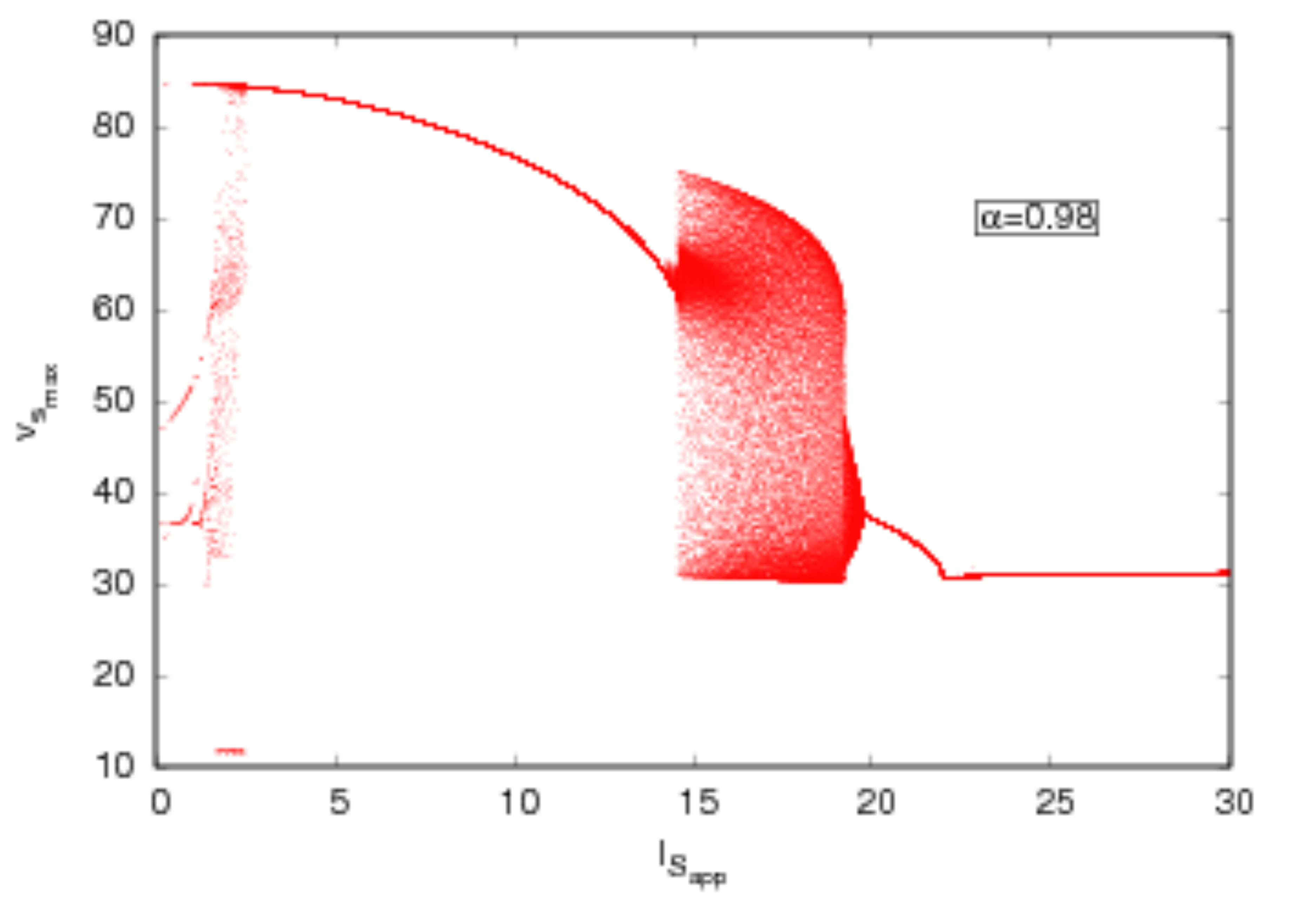}\\
\includegraphics[width=7cm,height=6cm]{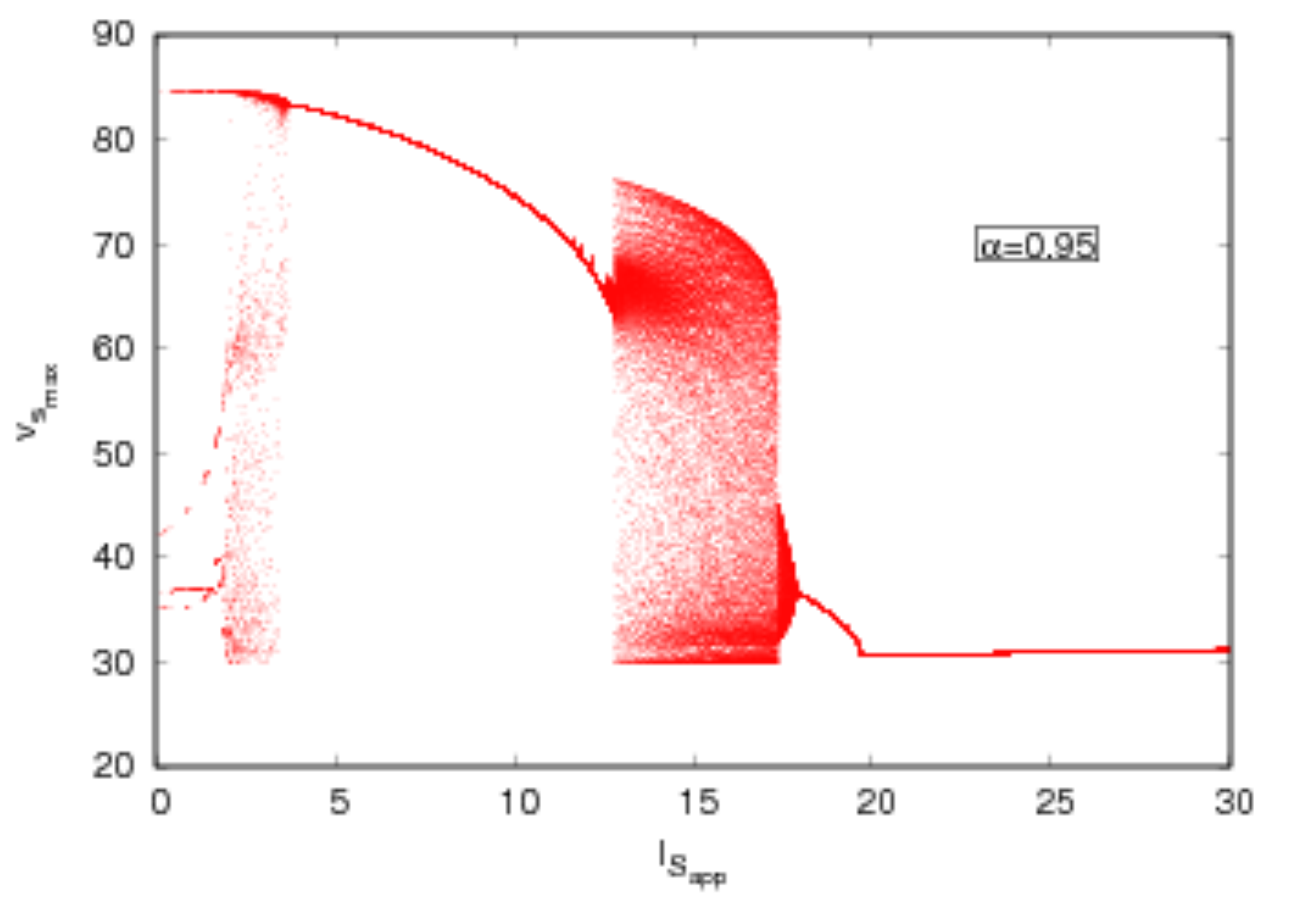}&
\includegraphics[width=7cm,height=6cm]{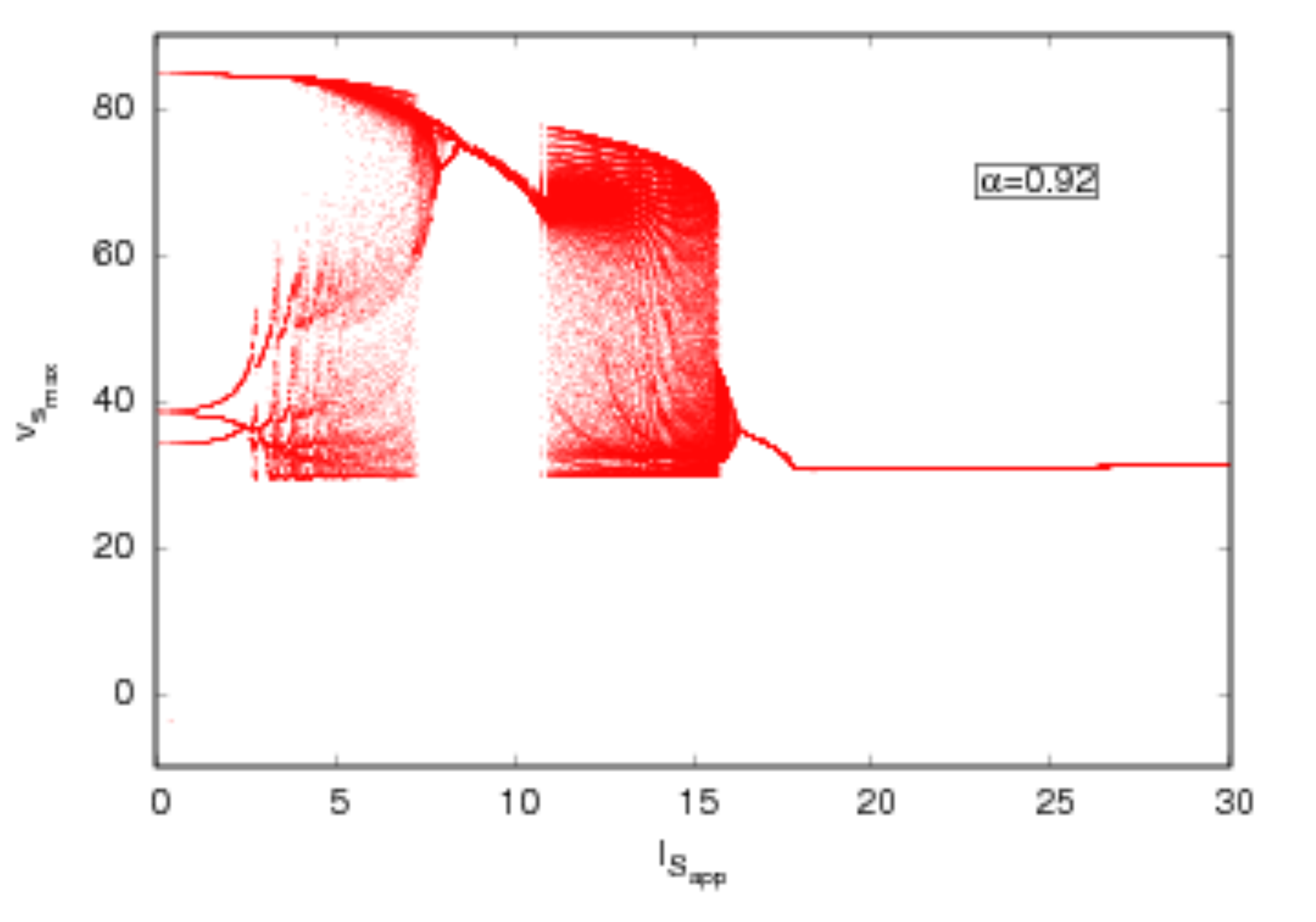}\\
\includegraphics[width=7cm,height=6cm]{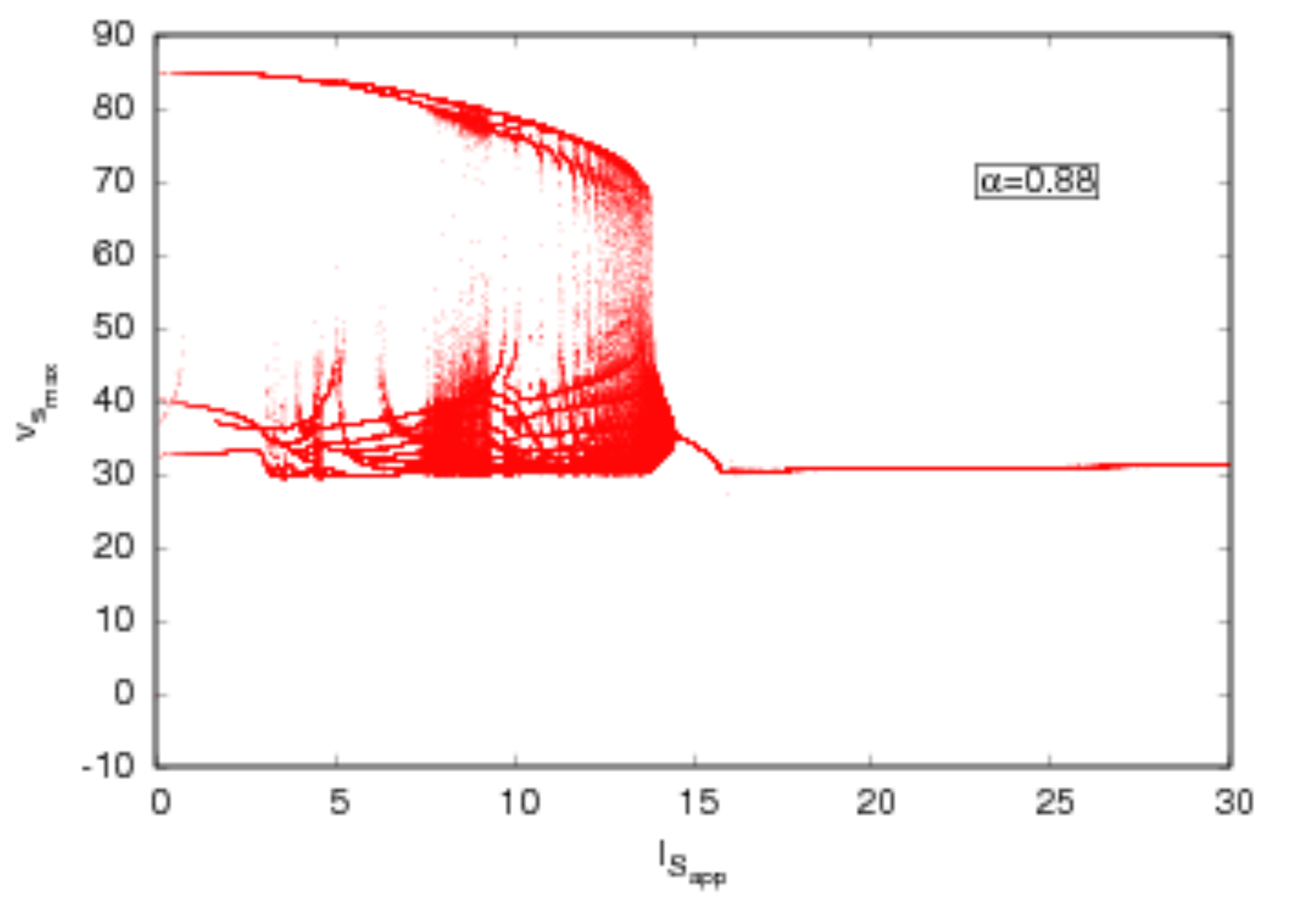}&
\includegraphics[width=7cm,height=6cm]{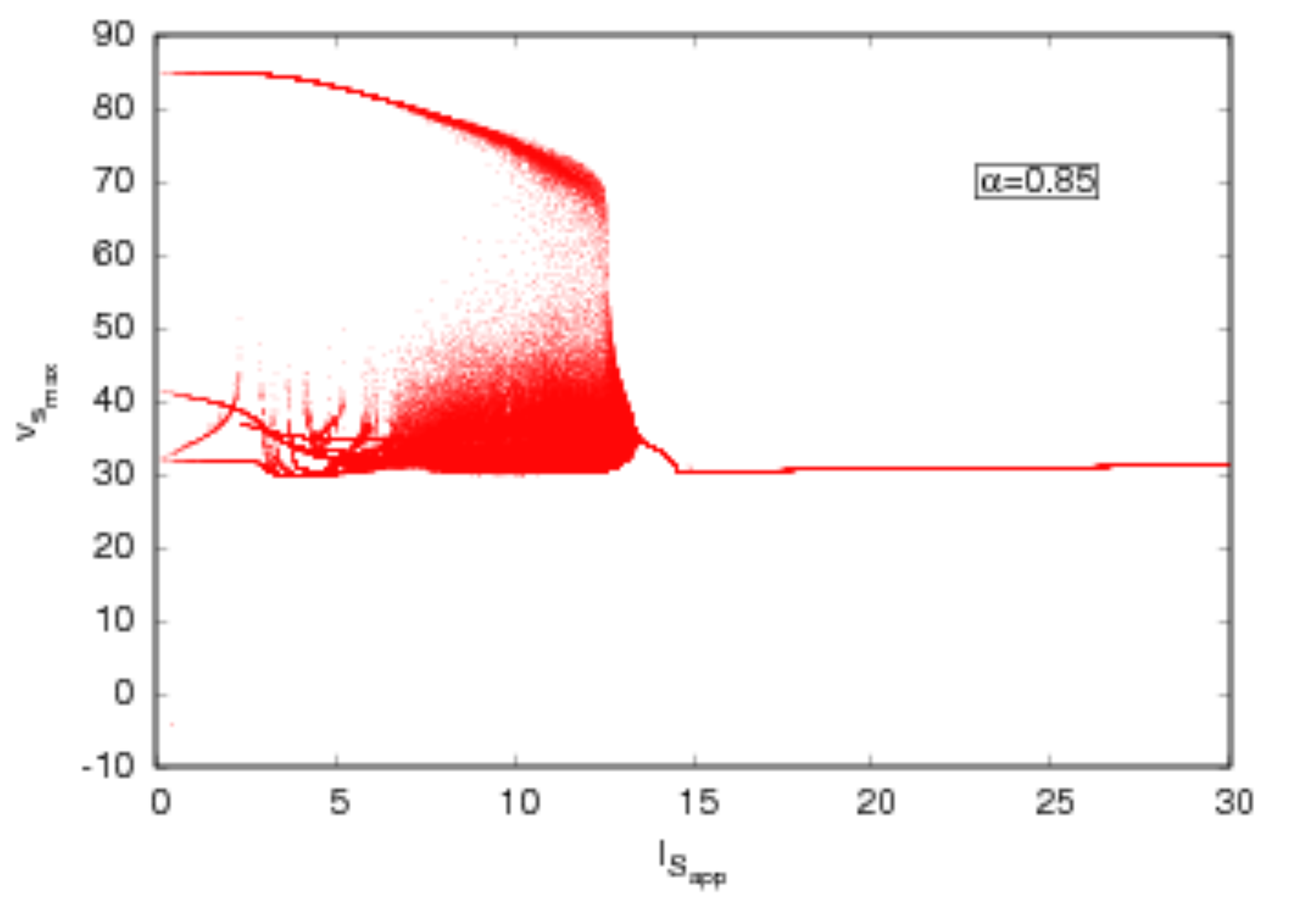}
\end{tabular} 
\caption{Bifurcation diagram when the parameter $\si{I_{S_{app}}}$  varied with the different values of the order $\alpha$. } \label{fig11}
\end{figure}
Moreover, this fact is true when we consider the $\si{I_{D_{app}}}$ as  the bifurcation parameter (Fig. \ref{fig12}).

\begin{figure}[!tbp] 
\centering
\begin{tabular}{cc}
\includegraphics[width=7cm,height=6cm]{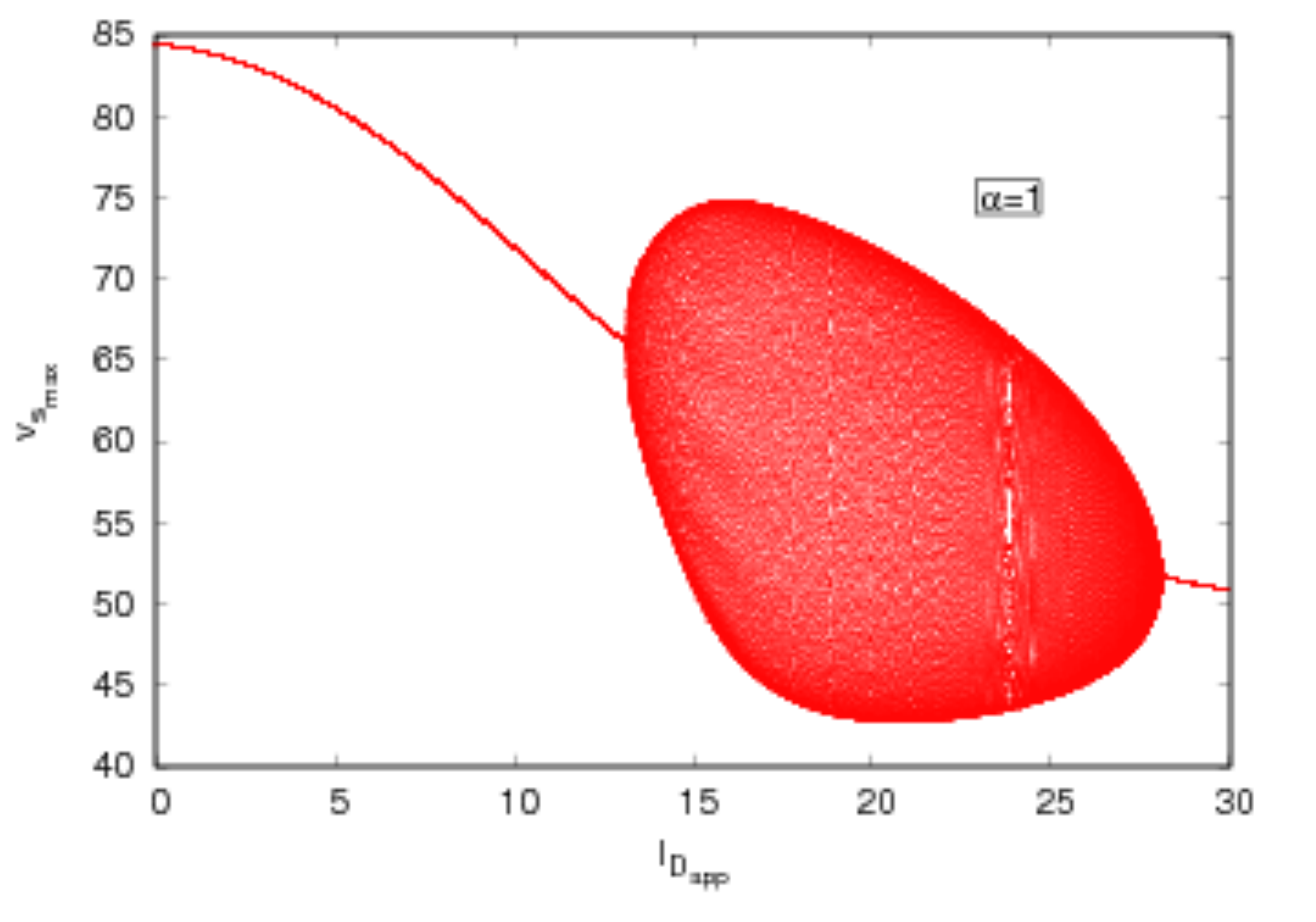}&
\includegraphics[width=7cm,height=6cm]{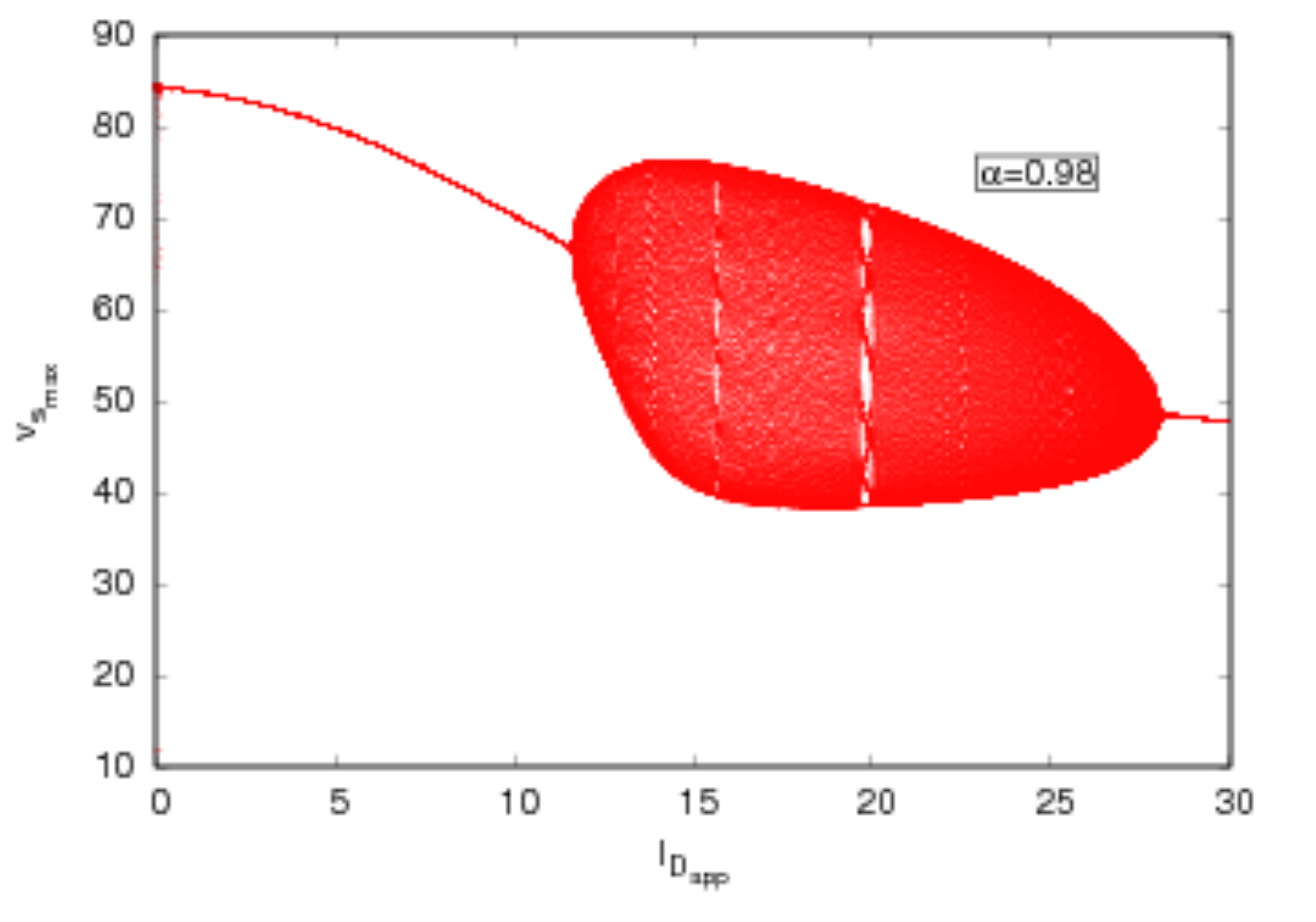}\\
\includegraphics[width=7cm,height=6cm]{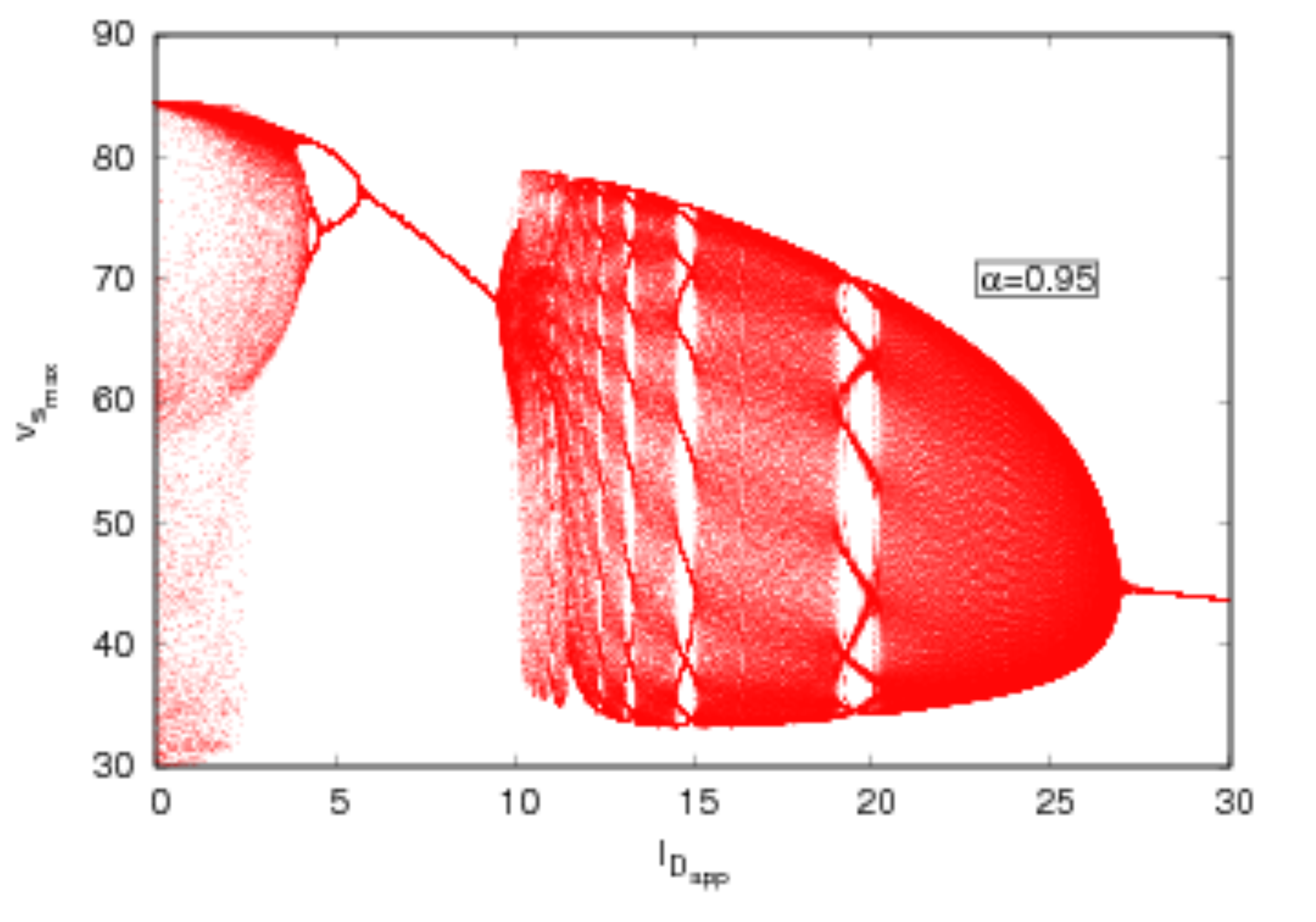}&
\includegraphics[width=7cm,height=6cm]{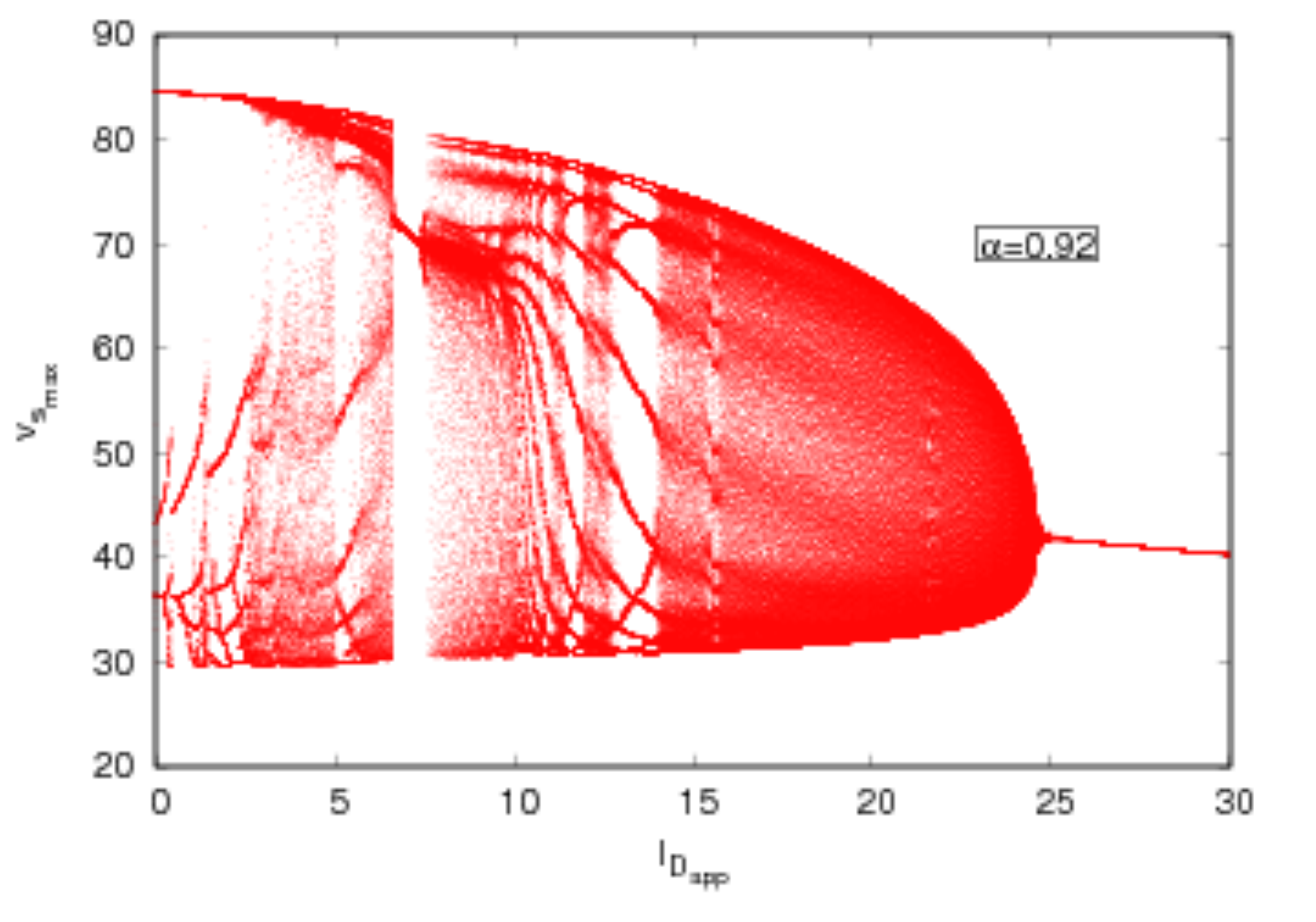}\\
\includegraphics[width=7cm,height=6cm]{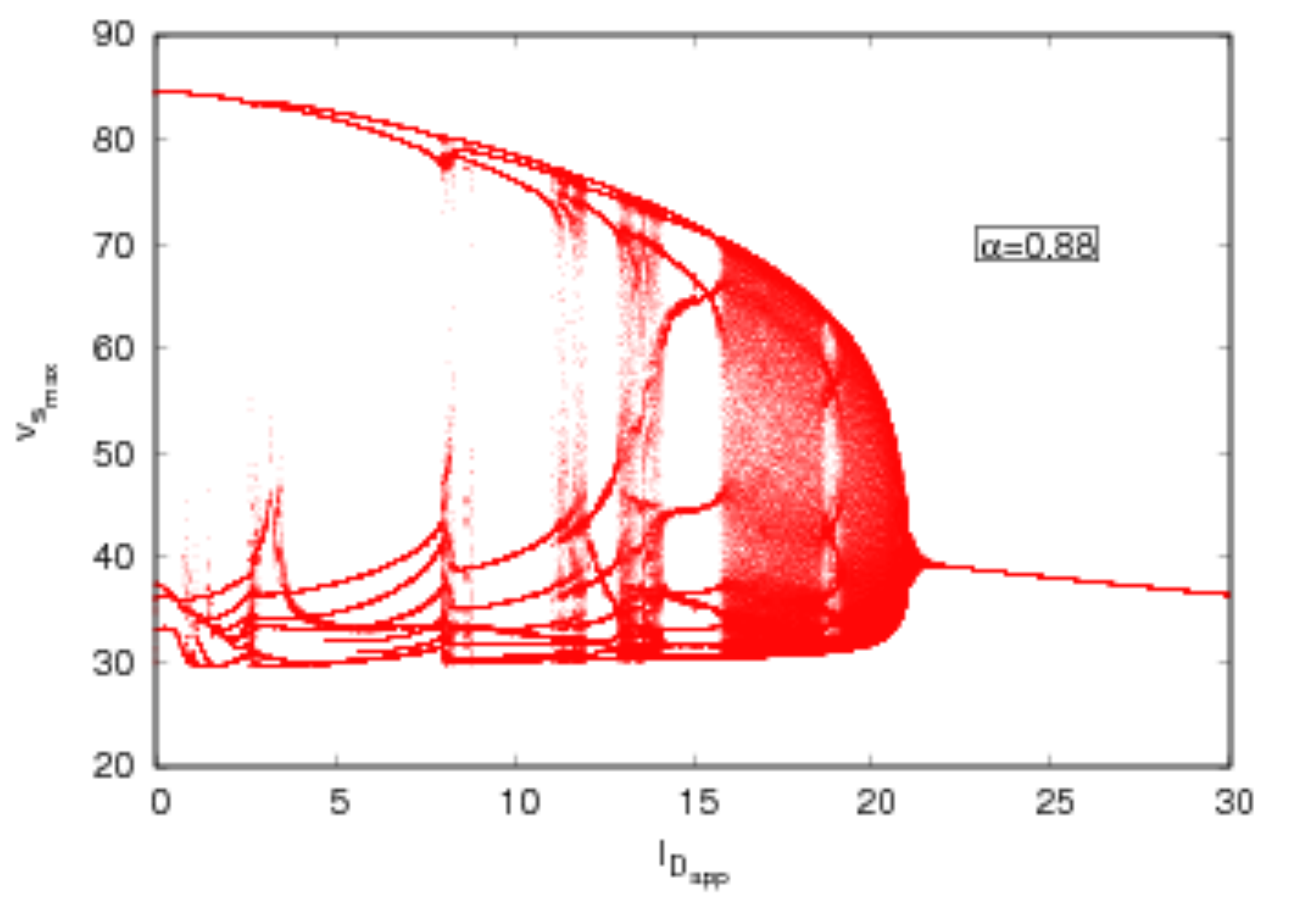}&
\includegraphics[width=7cm,height=6cm]{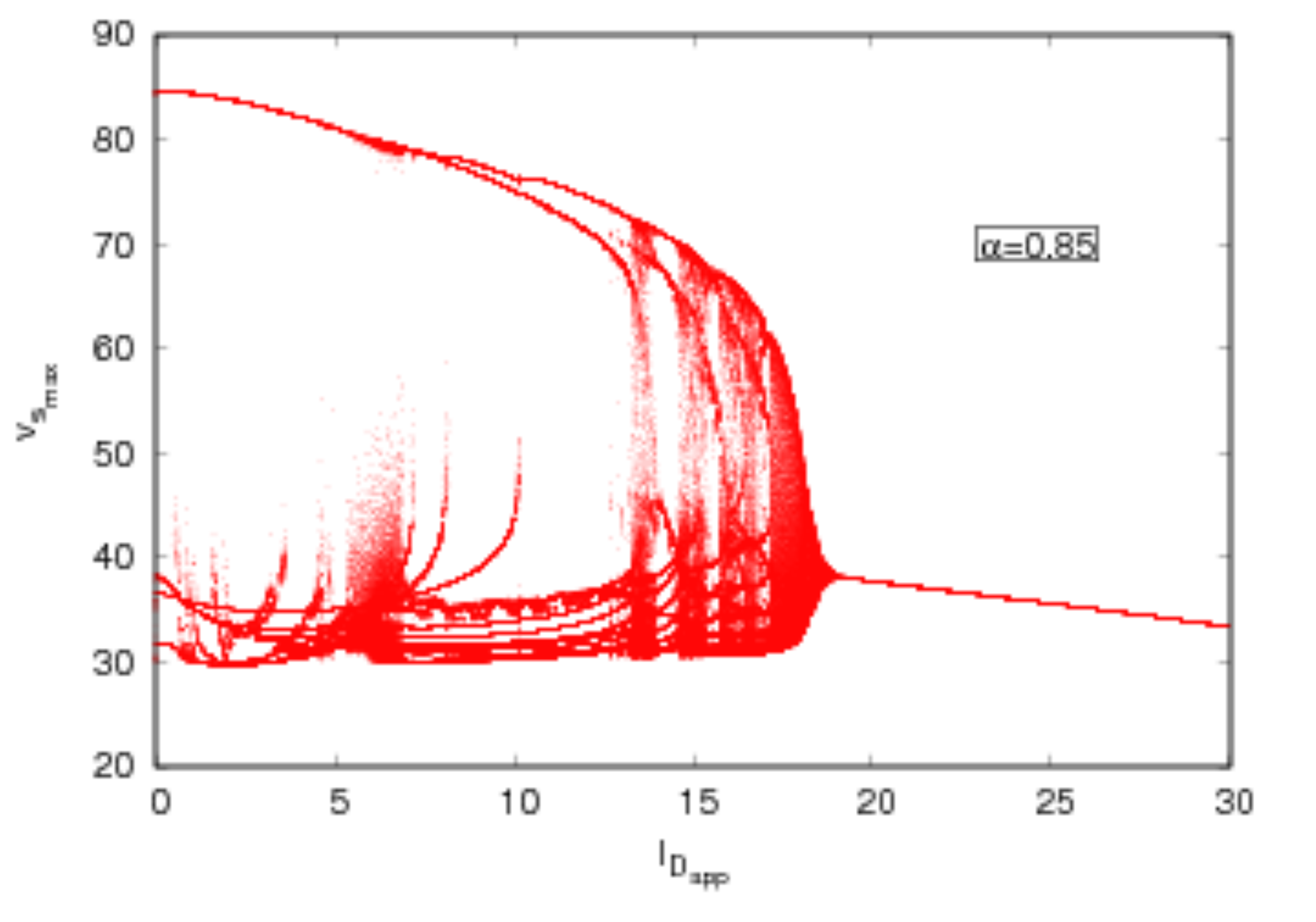}
\end{tabular} 
\caption{Bifurcation diagram when the parameter $\si{I_{D_{app}}}$  varied with the different values of the order $\alpha$. } \label{fig12}
\end{figure}

In Figs. \ref{fig11} and \ref{fig12}, two distinct regions are still visible, and we can clearly see chaotic as well as non-chaotic regions. The points where the regions change indicate fractional Hopf bifurcations (a loss of the equilibrium's stability).
As seen in Figs. \ref{fig11} and \ref{fig12}, the chaotic region starts to fade by reducing the value of $\alpha$  from 1 to 0.8. In addition, Figs. \ref{fig11}, \ref{fig12}, and \ref{fignew1} do not show any self-similarity.

\subsection{Stability analysis}\label{sub1}
In the following, we focus on the stability of the suggested fractional-order systems under specific conditions. Finding an analytical solution for the proposed model seems impossible in general. In the literature, investigating the stability for such large dynamics has always been challenging. However, for our analysis, we take into account a numerical scheme. We study the stability of the system based on a theorem for the Jacobian matrix of the system. Because of the non-linear nature of the system, we then obtain the equilibrium points numerically. 

\begin{defn}
The autonomous nonlinear differential system in the sense of Caputo is given as follows
\begin{equation}\label{e22}
{}^CD^{\alpha}x(t)=f(x(t)),
\end{equation}
where $x(t)\in \mathbb{R}^n$, $f(x(t))\in\mathbb{R}^{n\times 1}$ is a nonlinear function vector, $0< \alpha< 1$.
\end{defn}
\begin{defn}
The point $ \varepsilon^* $ is an equilibrium point of
system \eqref{e22}, if and only if $ f (\varepsilon^*) = 0$.
\end{defn}
\begin{defn}
The equilibrium point $ \varepsilon^* $ of \eqref{e22} is said to be

(1) stable if $ \forall \epsilon> 0$, there exists a $ \delta>0 $
such that $ \Vert x(t) - \varepsilon^*  \Vert< \epsilon $ holds for $ \forall x_0\in\lbrace z:\Vert z-\varepsilon^*\Vert<\delta\rbrace $ and $ \forall t > 0 $;

(2) asymptotically stable if the equilibrium point is locally stable and $ \lim_{t\rightarrow +\infty} x(t) = \varepsilon^*  $.
\end{defn}
\begin{thm}\label{rem2}
Consider the non-linear fractional-order autonomous system \eqref{e22}. Let  $ \varepsilon^* $ is an equilibrium point of the nonlinear system \eqref{e22}, then  equilibrium point $ \varepsilon^* $ is asymptotically stable if eigenvalues of Jacobian matrix $\lambda(J(\varepsilon^*))$ of the system satisfy 
$\left| \mathrm{arg} (\mathrm{spec}(J(\varepsilon^*)))\right| >\dfrac{\alpha \pi}{2}.$
Where  $\mathrm{spec}(J(\varepsilon^*))$ is the set of all eigenvalues of $J(\varepsilon^*)$\cite{kaslik2012nonlinear}. 
\end{thm}\label{rem1}

Since we cannot establish the stability conditions for this system in a general form, we attempt to find a numerical technique for certain values. One can consider any other values based on their needs. The related codes for the implementation of the methods can be seen in Appendix.

We use Matlab software and the $fsolve$ command to obtain the equilibrium points of the system. Based on Theorem (\ref{rem2}) we also compute the Jacobian matrix of the system to study the stability of the equilibrium points. Deu to chaotic behaviour of the system in the interval [0.93,1] we choose $\alpha=0.95$ and try to find some domain of ${I_{S_{app}}}$ and ${I_{D_{app}}}$ where the system \eqref{bif1} would be asymptotically stable.
Here, we evaluate two cases. First, we consider  $\alpha=0.95\in[0.93,1], $   $\si{I_{S_{app}}}=0.75,~  \si{I_{D_{app}}}=0,$ and $ g_c=2.1$.  
The equilibrium point of the system is
 \[\varepsilon^*_1\approx [1.0225, 0.9267, 0.9950, 0.0088 , 0.0149, 0.0117, 0.0456, 0.5355],\]
 and the eigenvalues of the Jacobian matrix are as 
 \[
 \mathrm{spec}(J(\varepsilon^*_1))\approx[ -3.3066, -2.4601, -1.0242,  -0.3830, -0.3123, -0.0749, 0.0007 + 0.0005i, 0.0007 - 0.0005i],
 \]
 For the second case, we consider
 $\si{I_{S_{app}}}=2.5,~  \si{I_{D_{app}}}=0,$ and   $g_c=2.1$. 
 in the case2, the equilibrium point is given by
 \[\varepsilon^*_2\approx [2.6161, 2.4916, 0.9929, 0.0239, 0.0170, 0.0132, 0.0448, 0.6876],\]
 and the eigenvalues  are as
 \[
 \mathrm{spec}(J(\varepsilon^*_2))\approx[ -3.2804, -2.3217, -0.9954,  -0.3453, -0.2856, -0.0748, 0.0039,  -0.0000]
 \]
 In both cases, $\lvert arg(\lambda _i)\rvert$  for $i=1,\cdots,8$, do not satisfying the conditions of the Theorem (\ref{rem2}) and then the equilibria $\varepsilon^*_1$ and $\varepsilon^*_2$ are not asymptotically stable.

\begin{figure}[!tbp] 
\centering
\begin{tabular}{cc} 
\includegraphics[width=7cm,height=6cm]{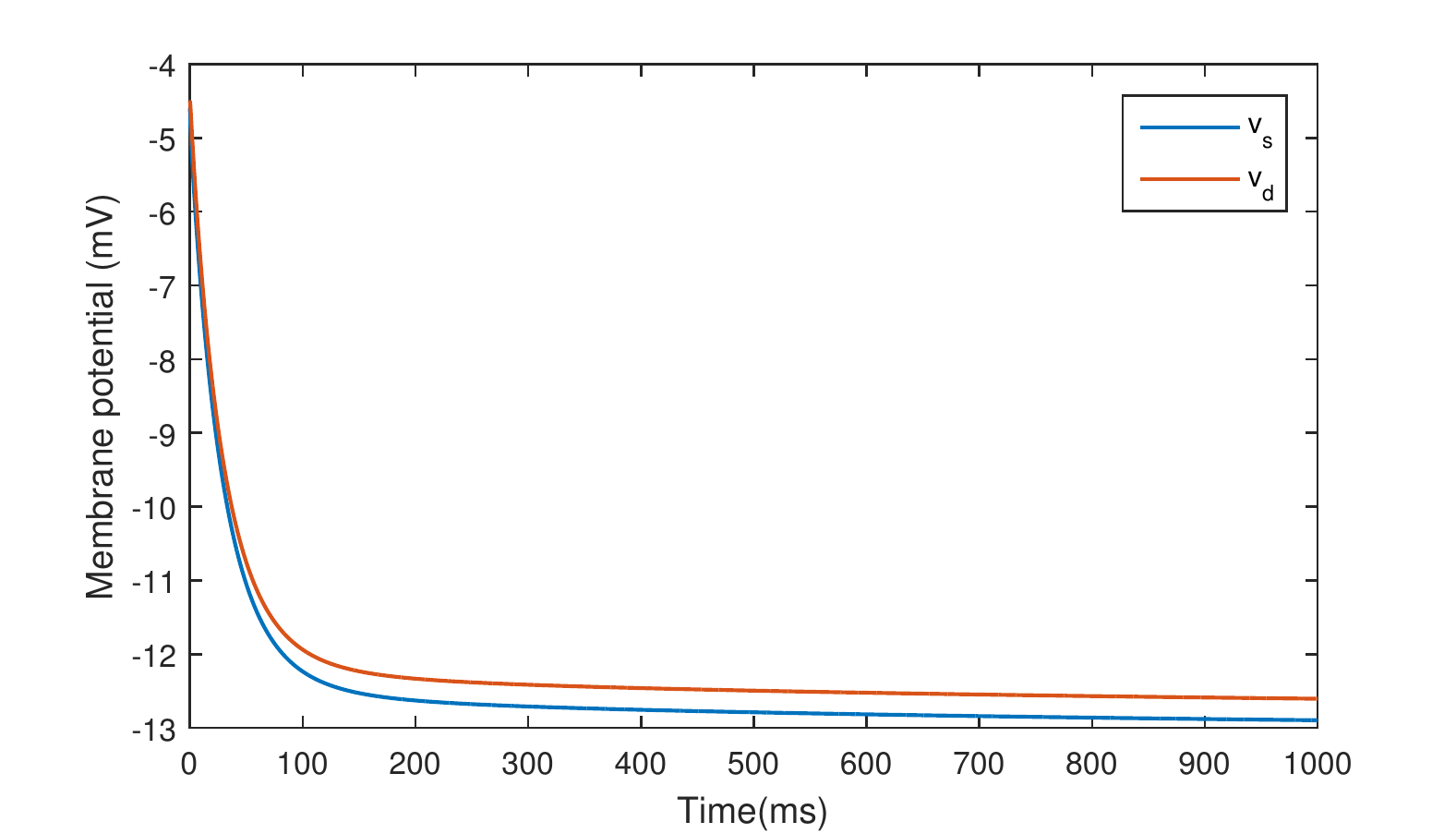} 
\includegraphics[width=7cm,height=6cm]{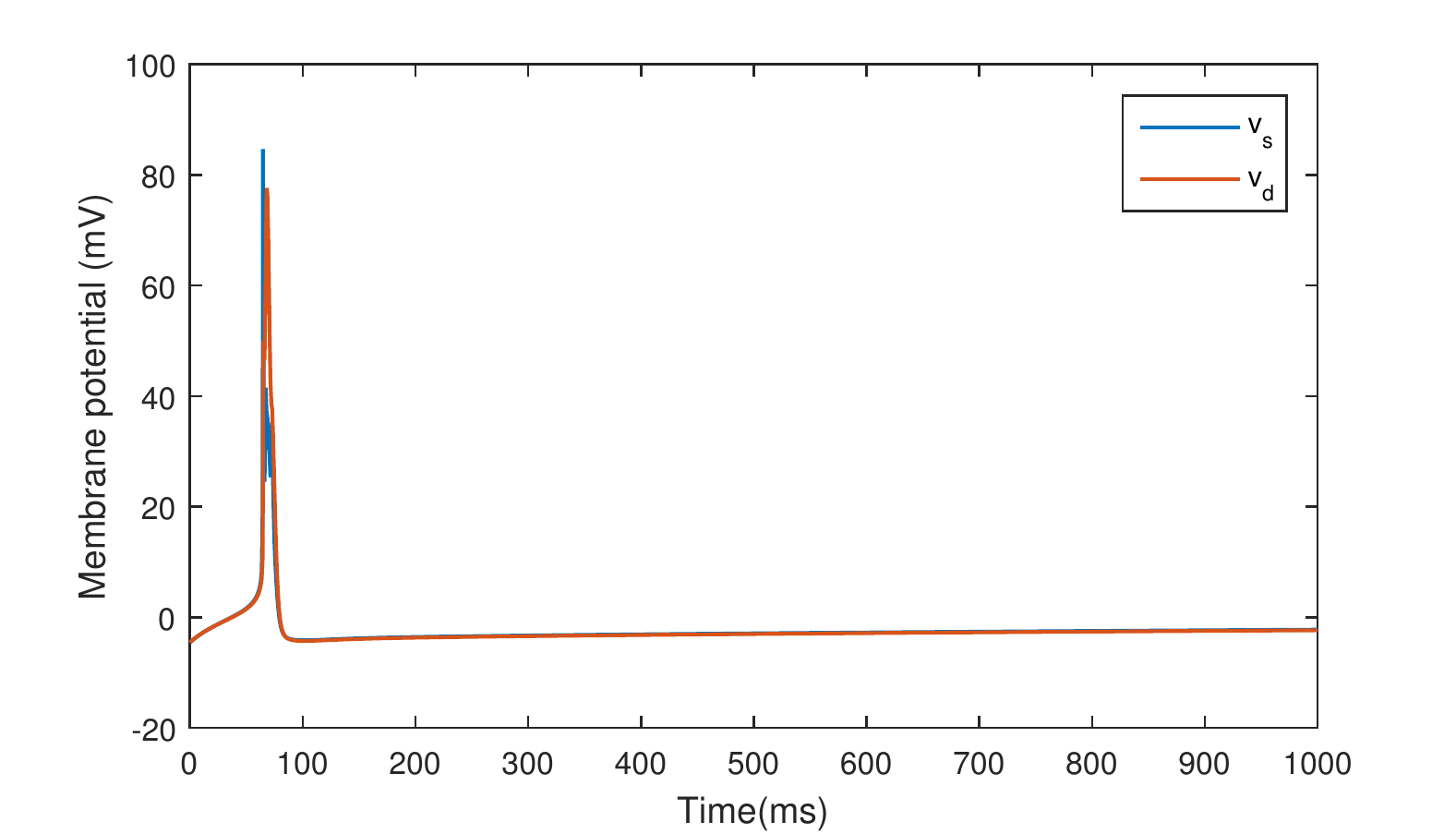} 
\end{tabular}
\caption{ Memberance potential for $\alpha=0.95$, $\si{I_{D_{app}}}=0$,  $\si{I_{S_{app}}}=-1.2579$ (left), $\si{I_{S_{app}}}=0.0268$ (right).}\label{fig13} 
\end{figure}

\begin{figure}[!tbp] 
\centering
\begin{tabular}{cc} 
\includegraphics[width=7cm,height=6cm]{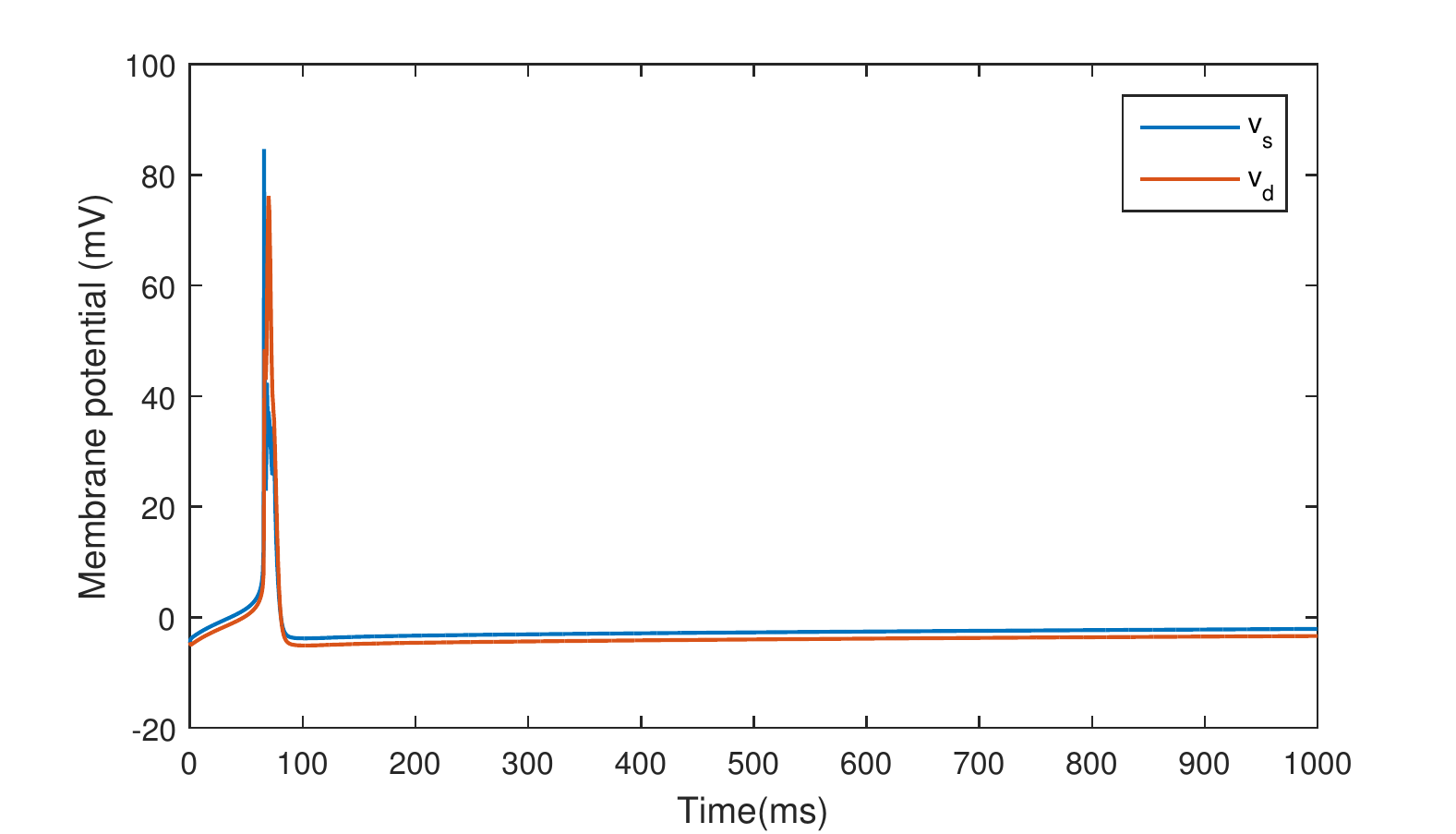} 
\includegraphics[width=7cm,height=6cm]{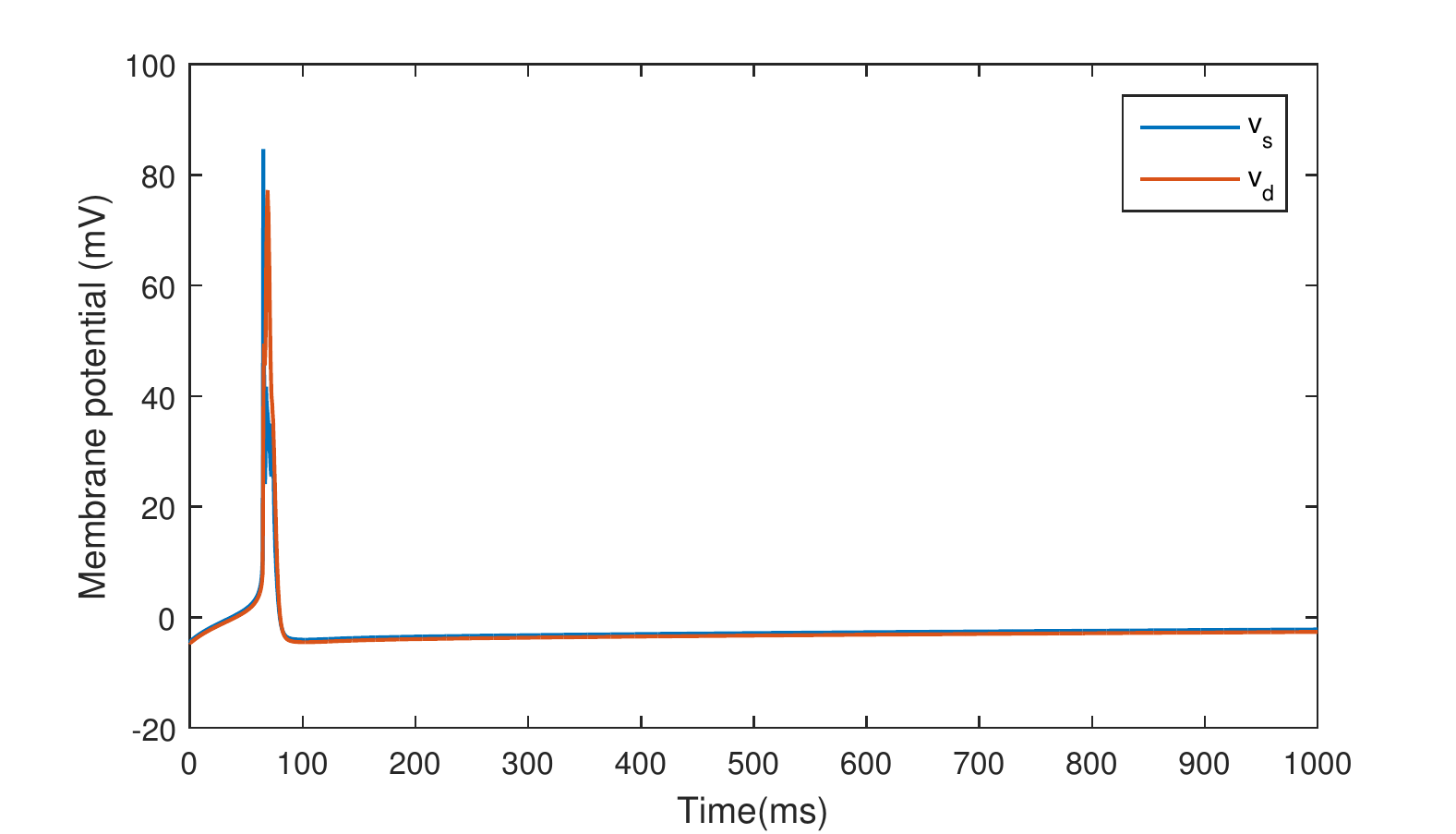}
\end{tabular} 
\caption{ Memberance potential for $\alpha=0.95$, $\si{I_{S_{app}}}=2.5$ and $\si{I_{D_{app}}}=-2.5471$ (left), $\si{I_{S_{app}}}=0.75$ and $\si{I_{D_{app}}}=-0.7449$ (right).}\label{fig15} 
\end{figure}

Now we try to find some domains of $ \si{I_{S_{app}}}$ and $ \si{I_{D_{app}}}$  near the standard value of $ \si{I_{S_{app}}}=2.5 $ and $ \si{I_{D_{app}}}=0$ to make the system \eqref{bif1} asymptotically stable.
At first, we consider $\si{I_{S_{app}}}$ as a parameter that changes in the interval $[-4, 4]$. Using Remark(\ref{rem2}) and an increment equal to $0.0001$ in $\si{I_{S_{app}}}$ we found the interval $[-1.2579, 0.0268]$ where the system does not show chaotic behaviour and the equilibria $\varepsilon^*$ is asymptoticalliy stable (see Figs.\ref{fig13} ).
Secondely, taking $\si{I_{D_{app}}}$ as a parameter changing in the interval $[-4, 4]$ with $\si{I_{S_{app}}}=2.5 $ and $ \si{g_c}=2.1 $, using the Remark(\ref{rem2}) and an increment equal to $0.0001$ in $\si{I_{D_{app}}}$ we found the interval $[-4, -2.5471]$ where the system does not show chaotic behaviour and also the equilibria $\varepsilon^*$ of the system is asmptoticalliy stable (see Fig.\ref{fig15}).
Finally, we let $\si{I_{D_{app}}}$ as parameter in the interval [-4, 4] with $\si{I_{S_{app}}}=0.75 $ and $ \si{g_c}=2.1 $. Using the Remark(\ref{rem2}) and an increment equal to $0.0001$ in $\si{I_{D_{app}}}$ we found the interval $[-4, -0.7449]$ where the system doesn't show chaotic behaviour and also the equilibria $\varepsilon^*$ is asmptoticalliy stable (see Fig.\ref{fig15}).

\section{Discussion}
This work investigated the bifurcation analysis of a fractional-order model of a two-compartment  CA3  hippocampal pyramidal cell. Chaotic regions were achieved for different values of the fractional derivative order and different injection currents. 
Since the membrane capacitance has been considered ideal, the fractional-order model has deemed as a more accurate description of physical processes underlying a long-range memory behavior. We can model the dielectric losses in the membrane with fractional capacitance and observe a change in the relative refractory period, as shown in  Fig. \ref{figtrpcur}. In fact, Fig. \ref{figtrpcur}  shows a comparison between fractional-order and integer-order cases. In fractional case  $\alpha=0.85$ the refractory period is smaller than the integer case $\alpha=1$. Moreover, based on Figs. \ref{fig11} and \ref{fig12}, it can be concluded that the new fractional-order system's chaotic region collapses when the derivative order approaches from 1 to 0. The obtained results can be considered as help to control relevant diseases
caused by maximal injection currents abnormality. 
\renewcommand{\thefigure}{S\arabic{figure}}
\setcounter{figure}{0}
\section*{Appendix}
\begin{figure}
\centering
\includegraphics[width=\textwidth]{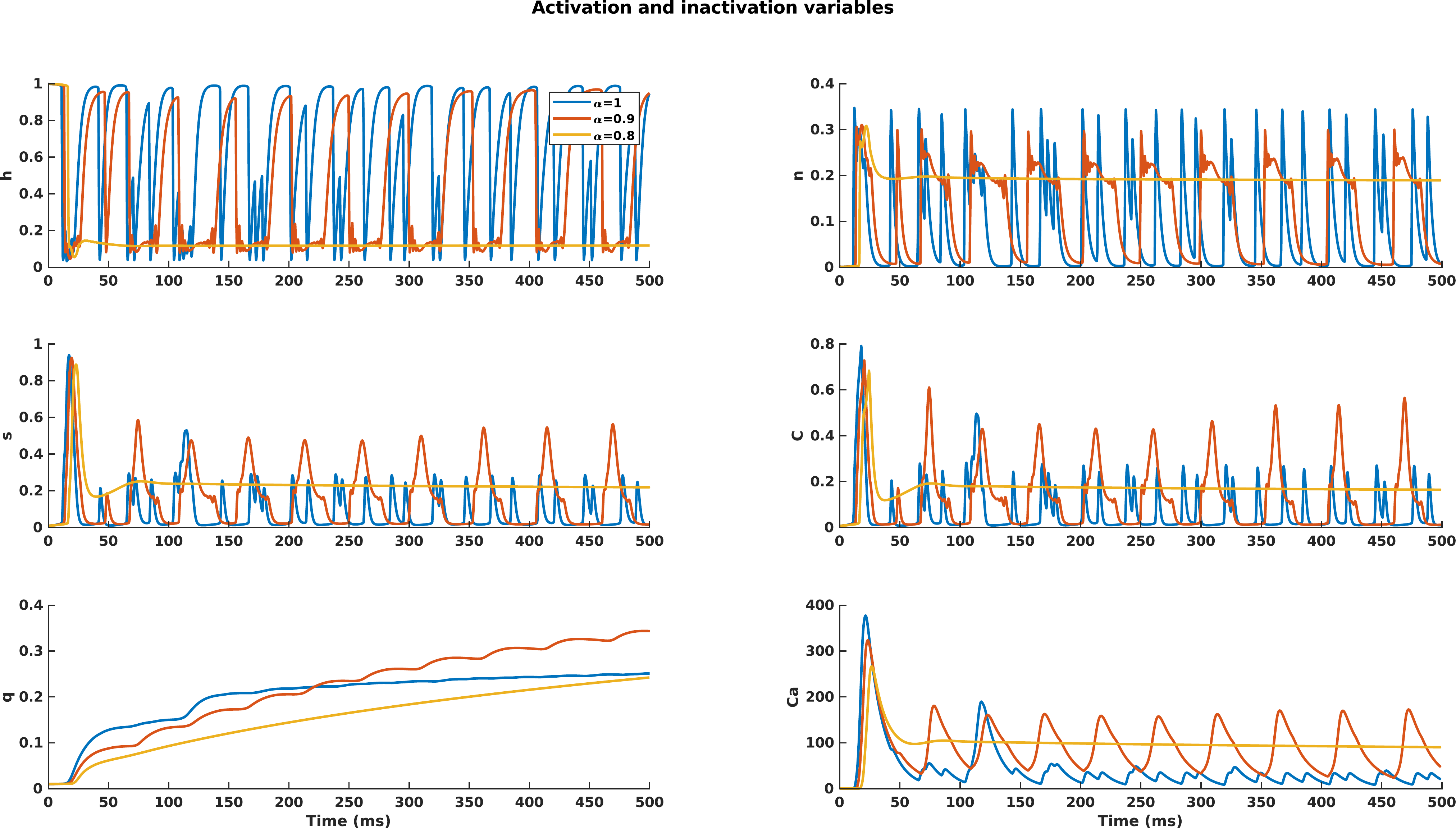}
\caption{Activation and inactivation variables when $\si{I_{S_{app}}}=2.5$ and $\si{I_{D_{app}}}=0$.}\label{fig25var}
\end{figure}

\begin{figure}
\centering
\includegraphics[width=\textwidth]{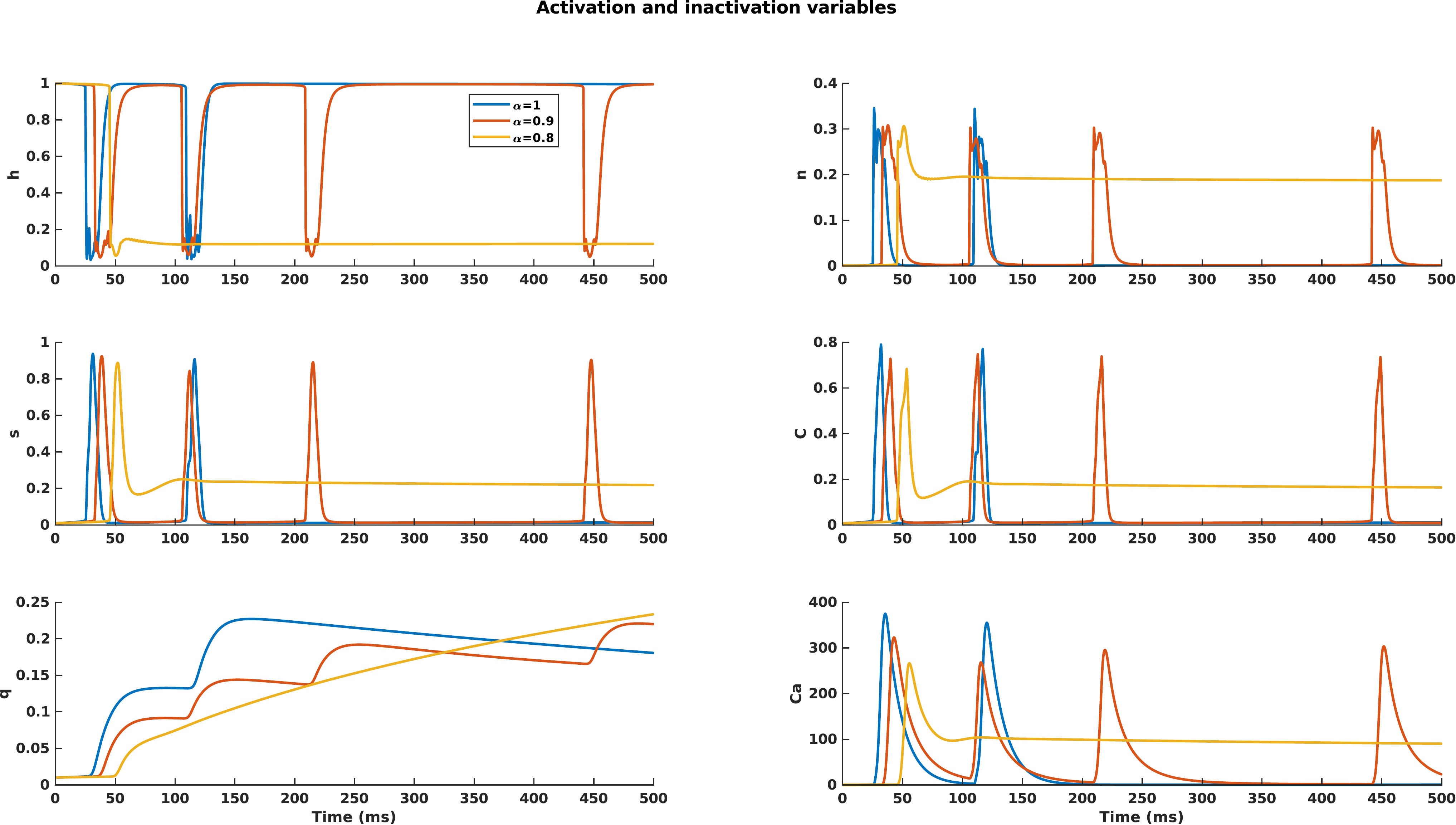}
\caption{Activation and inactivation variables when $\si{I_{S_{app}}}=0.75$ and $\si{I_{D_{app}}}=0$.}\label{fig075var}
\end{figure}

\begin{figure}
\centering
\includegraphics[width=\textwidth]{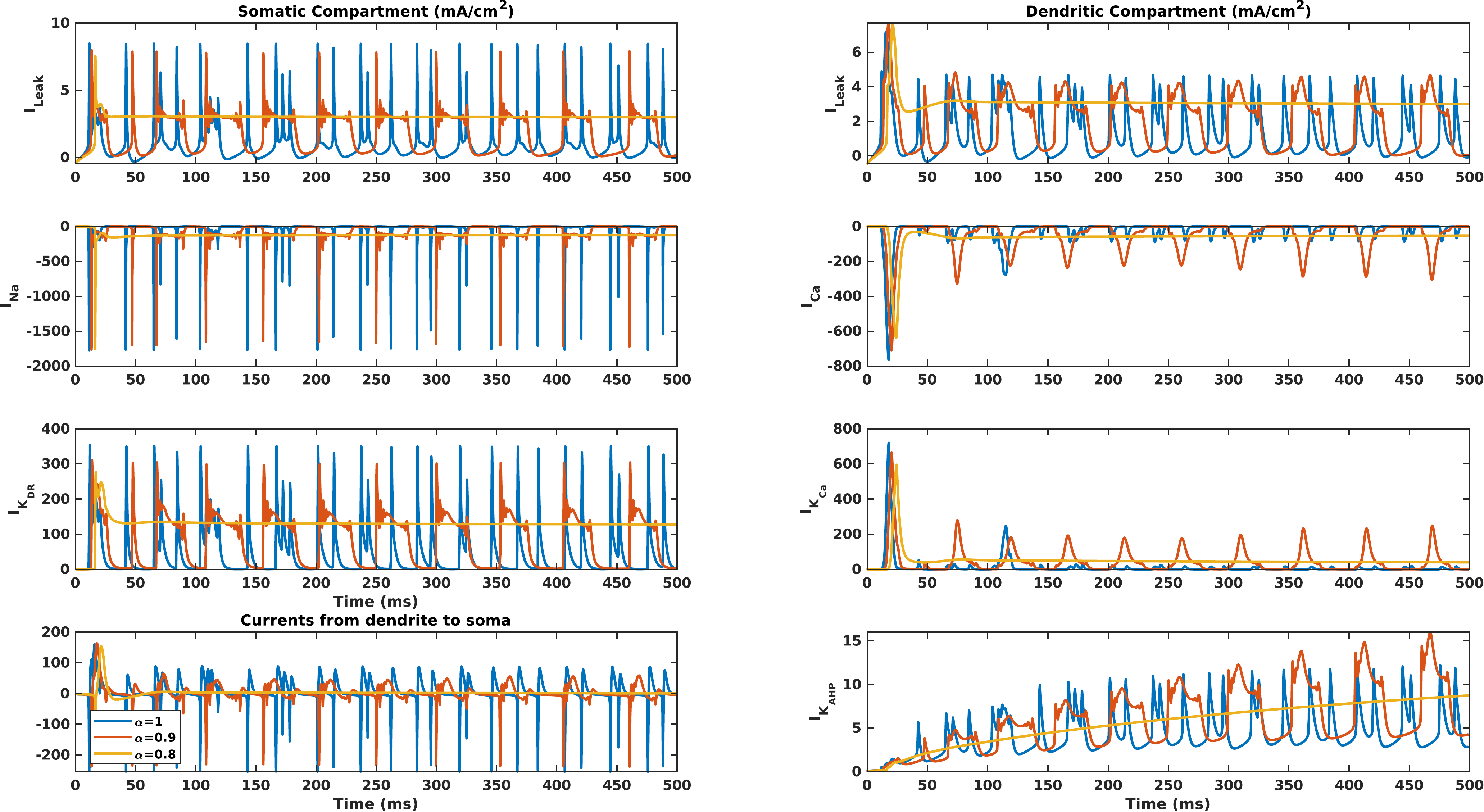}
\caption{Time course of the currents in the somatic compartment, from somatic to dendritic compartment, and in the dendritic compartment when $\si{I_{S_{app}}}=2.5$ and $\si{I_{D_{app}}}=0$..}\label{fig25cur}
\end{figure}

\begin{figure}
\centering
\includegraphics[width=\textwidth]{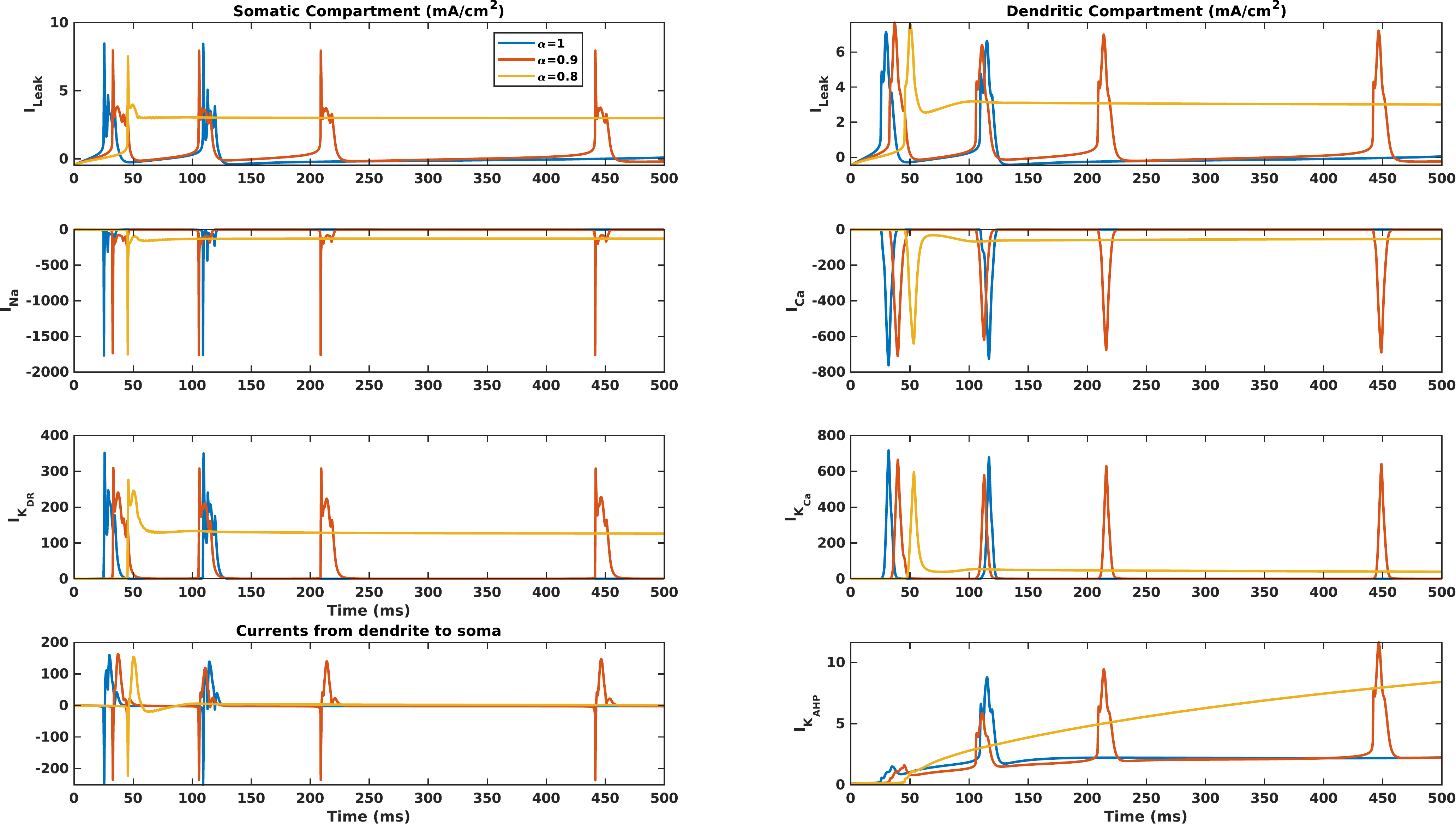}
\caption{Time course of the currents in the somatic compartment, from somatic to dendritic compartment, and in the dendritic compartment when $\si{I_{S_{app}}}=0.75$ and $\si{I_{D_{app}}}=0$.}\label{fig075cur}
\end{figure}

\begin{figure}
\centering
\includegraphics[width=\textwidth]{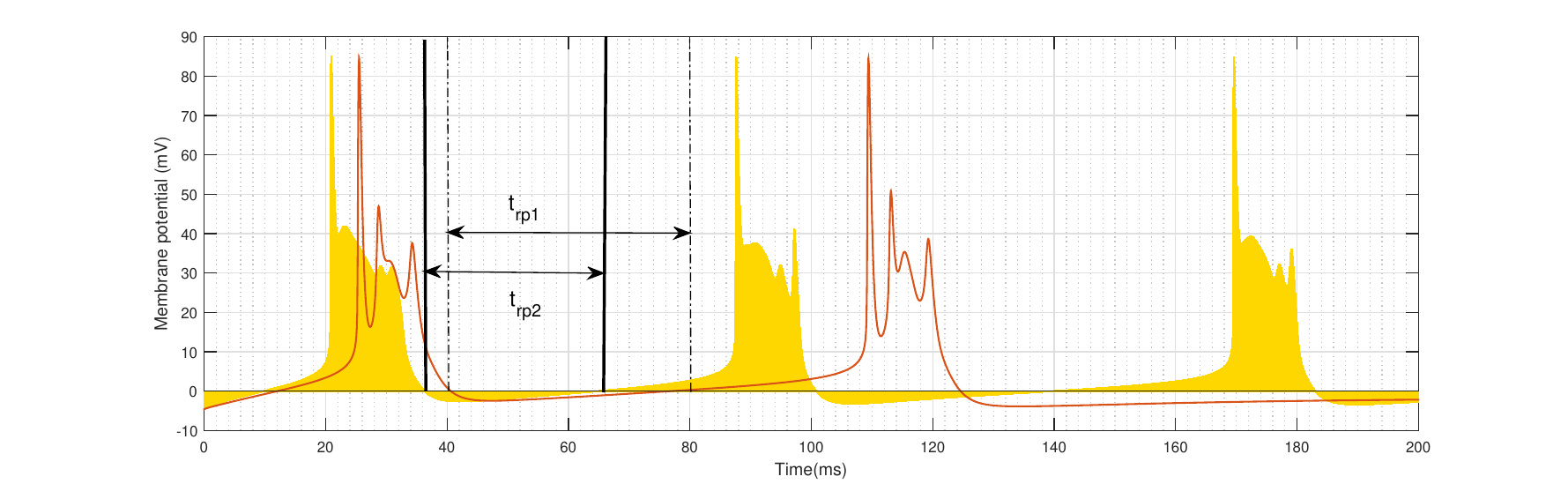}
\caption{refractory period ($t_{rp2}$) for $\alpha=0.85$ (yellow) and refractory period ($t_{rp1}$) for $\alpha=1$ (red).}\label{figtrpcur}
\end{figure}

In this section, we list the notations and the value of the parameters used for the numerical results (Table \ref{t2}), and illustrate the supplementary figures (Figs \ref{fig25var}-\ref{fig075cur}).
\section{Matlab code }
\begin{lstlisting}[language=octave]
clc 
close all
%  y=[Vs, Vd, h, n, s, c, q, Ca]
% Id=I_{D_{app}} and Is= I_{S_{app}}, Ild,Ils=I_{Leak}

y0= [-4.6 ;-4.5 ;0.999 ;0.001 ;0.009 ;0.007; 0.01; 0.2]
Id = 0:.01:4;
i=0
alpha=0.95
for ct = 1:numel(Id)
   y(:,ct)=fsolve(@(y)myfun(y,Id(ct)),y0);
   df=NumJacob(@(y)myfun(y,Id(ct)),y(:,ct));
  EIG= eig(df);
  angle(EIG);
 idx = find(abs(angle(EIG))<=(alpha*pi/2));
    if (numel(idx~=0))
     disp(['is not asymptotically stable']);
       else
     i=i+1
     disp(['is asymptotically stable'])
  end
end
\end{lstlisting}
\begin{lstlisting}[language=octave]
function f=test_func(y,Id)
gL = .1; gNa = 30; gKDR =15; gCa = 10; gKAHP = 0.8; gKC =15; 
VNa = 120; VCa = 140; VK = -15; VL = 0;
Isyn = 0; p = 0.5; alpha=.95; Rm=10; tm=30; Cm=tm^alpha/Rm;
Is=2.5; gc=2.1; 
Ils = gL*(y(1)-VL);
minf = am(y(1))./(am(y(1))+bm(y(1)));
INa = gNa*minf.^2.*y(3).*(y(1)-VNa);
IKDR = gKDR*y(4)*(y(1)-VK);
Ild = gL*(y(2)-VL);
ICa = gCa*y(5).^2.*(y(2)-VCa);
K=(1.073*sin(0.003453*y(8)+0.08095)+ 0.08408*sin(0.01634*y(8)-2.34)+ 0.01811*sin(0.0348*y(8)-0.9918));
IKC = gKC*y(6)*K*(y(2)-VK);
IKAHP = gKAHP*y(7).*(y(2)-VK);

f = [(-Ils - INa - IKDR + (gc/p)*(y(2)-y(1)) + Is/p)/Cm;...
 (-Ild-ICa-IKAHP-IKC-Isyn/(1-p)+(gc/(1-p))*(y(1)-y(2))+Id/(1-p))/Cm;...
 ah(y(1))*(1-y(3)) - bh(y(1))*y(3);...
 an(y(1))*(1-y(4)) - bn(y(1))*y(4);...
 as(y(2))*(1-y(5)) - bs(y(2))*y(5);...
 (aq(y(8))-y(7))/bq(y(8));
  -0.13*ICa - 0.075*y(8)];
return 

function val = am(v)
v=v-60;
val = 0.32*(-46.9-v)/(exp((-46.9-v)/4)-1);

function val = bm(v)
v=v-60;
val = 0.28*(v+19.9)/(exp((v+19.9)/5)-1);

function val = an(v)
v=v-60;
val = 0.016*(-24.9-v)/(exp((-24.9-v)/5)-1);

function val = bn(v)
v=v-60;
val =0.25*exp(-1-0.025*v);

function val = ah(v)
v=v-60;
val = 0.128*exp((-43-v)/18);

function val = bh(v)
v=v-60;
val =4/(1+exp((-20-v)/5));

function val = as(v)
v=v-60;
val =1.6/(exp(-0.072*(v-5))+1);

function val = bs(v)
v=v-60;
val =0.02*(v+8.9)/(exp((v+8.9)/5)-1);
function val = aq(v)

val =(0.7894*exp(0.0002726*v))-(0.7292*exp(-0.01672*v));
function val = bq(v)

val =(657.9*exp(-0.02023*v))+(301.8*exp(-0.002381*v));
function val = ac(v)
v=v-60;
val =(1.0/(1+exp((-10.1-v)/0.1016)))^0.00925;

function val = bc(v)
v=v-60;
val =3.627*exp(0.03704*v);

\end{lstlisting}
\begin{table}[th!]
\begin{center}
\caption{Notations of the model and maximal conductance parameters.}\label{t2}
\begin{tabular}{c c c}
\hline
Notations & Description  &  Values \\
\hline
\si{I_{Leak}}& Leak currents & ---\\
\si{I_{Na}}&Sodium current & ---\\
\si{I_{K_{DR}}}&Delayed rectifier K current &  ---\\
\si{I_{DS}}&Coupling current & ---\\
\si{I_{SD}}&Coupling current &  ---\\
\si{I_{S_{app}}}& Somatic current injection &  ---\\
\si{I_{Ca} }& Dendritic calcium current &  ---\\
\si{I_{K_{Ca}}}& Voltage and calcium dependent potassium current & ---\\
\si{I_{K_{AHP}}}& Ca dependent AHPK current& ---\\
\si{I_{D_{app}}}& Dendritic current injection& ---\\
\si{V_s}& Somatic membrance potentials & ---\\
\si{V_d}&Dendritic membrance potentials&---\\
\si{V_{Na}}&Reversal potentials, relative to resting potential, for$Na^+$ &60 \\
\si{V_{K}}&Reversal potentials, relative to resting potential, for $K^+$&-75\\
\si{V_{Ca}}&Reversal potentials, relative to resting potential, for $Ca^+$&80 \\
\si{V_{L}}& Reversal potentials, relative to resting potential, related to $I_{L}$& -60\\
\si {p }&The size of the axosomatic compartment as a proportion of the entire cell& 0.5  \\
\si {1-p }&The size of the dendritic compartment as a proportion of the entire cell& 0.5 \\
\si{C_m}& The capacitance & 3  \\
\si{g_{Na}}&Maximum conductance for Na &30 \\
\si{g_{K_{DR}}}&Maximum conductance for delayed rectifier K current& 15 \\
\si{g_{K_{Ca}}}& Maximum conductance for voltage and Ca dependent K current (fast)&  15 \\
\si{g_{K_{AHP}}}& Maximum conductance for Ca dependent AHPK current (slow)&   0.8\\
\si{g_{Ca}}&Maximum conductance for Ca&  10  \\
\si{g_{L}}& Leakage conductance &  0.1  \\
\si{g_{c} }&Strength of coupling &  2.1\\
\hline
\end{tabular}
\end{center}
\end{table}
\maketitle


\clearpage
\bibliographystyle{plain}
\bibliography{references.bib}

\begin{thebibliography}{10}

\bibitem{bib3}
Laura~A Atherton, Luke~Y Prince, and Krasimira Tsaneva-Atanasova.
\newblock Bifurcation analysis of a two-compartment hippocampal pyramidal cell
  model.
\newblock {\em Journal of computational neuroscience}, 41(1):91--106, 2016.

\bibitem{bib24a}
Donato Cafagna and Giuseppe Grassi.
\newblock Bifurcation and chaos in the fractional-order chen system via a
  time-domain approach.
\newblock {\em International Journal of Bifurcation and Chaos},
  18(07):1845--1863, 2008.

\bibitem{chen2008chaotic}
Juhn-Horng Chen and Wei-Ching Chen.
\newblock Chaotic dynamics of the fractionally damped van der pol equation.
\newblock {\em Chaos, Solitons \& Fractals}, 35(1):188--198, 2008.

\bibitem{chen2008nonlinear}
Wei-Ching Chen.
\newblock Nonlinear dynamics and chaos in a fractional-order financial system.
\newblock {\em Chaos, Solitons \& Fractals}, 36(5):1305--1314, 2008.

\bibitem{bib24}
Marius-F. Danca and Nikolay Kuznetsov.
\newblock Matlab code for lyapunov exponents of fractional-order systems.
\newblock {\em International Journal of Bifurcation and Chaos}, 28(05):1850067,
  2018.

\bibitem{de2014review}
Edmundo~Capelas de~Oliveira and Jos{\'e}~Ant{\'o}nio Tenreiro~Machado.
\newblock A review of definitions for fractional derivatives and integral.
\newblock {\em Math. Probl. Eng.}, 2014, 2014.

\bibitem{devaney1990chaos}
Robert~L Devaney.
\newblock {\em Chaos, fractals, and dynamics: computer experiments in
  mathematics}.
\newblock Addison-Wesley Longman Publishing Co., Inc., 1990.

\bibitem{dithelm}
Kai Diethelm and Neville~J Ford.
\newblock Analysis of fractional differential equations.
\newblock {\em Journal of Mathematical Analysis and Applications},
  265(2):229--248, 2002.

\bibitem{bib4}
G~Bard Ermentrout and David~H Terman.
\newblock {\em Mathematical foundations of neuroscience}, volume~35.
\newblock Springer Science \& Business Media, 2010.

\bibitem{bib5}
Gabbiani F and Cox S~J.
\newblock {\em Mathematics for Neuroscientists}.
\newblock Elsevier, 2017.

\bibitem{gao2005chaos}
Xin Gao and Juebang Yu.
\newblock Chaos in the fractional order periodically forced complex duffing’s
  oscillators.
\newblock {\em Chaos, Solitons \& Fractals}, 24(4):1097--1104, 2005.

\bibitem{bib28}
Roberto Garrappa.
\newblock Numerical solution of fractional differential equations: A survey and
  a software tutorial.
\newblock {\em Mathematics}, 6(2):16, 2018.

\bibitem{hartley1995chaos}
Tom~T Hartley, Carl~F Lorenzo, and H~Killory Qammer.
\newblock Chaos in a fractional order chua's system.
\newblock {\em IEEE Transactions on Circuits and Systems I: Fundamental Theory
  and Applications}, 42(8):485--490, 1995.

\bibitem{bib8}
Podlobny I.
\newblock {\em Fractional Differential Equations}.
\newblock Academic press, 1999.

\bibitem{kandel2000principles}
Eric~R Kandel, James~H Schwartz, Thomas~M Jessell, Steven Siegelbaum, A~James
  Hudspeth, and Sarah Mack.
\newblock {\em Principles of neural science}, volume~4.
\newblock McGraw-hill New York, 2000.

\bibitem{kaslik2012nonlinear}
Eva Kaslik and Seenith Sivasundaram.
\newblock Nonlinear dynamics and chaos in fractional-order neural networks.
\newblock {\em Neural Networks}, 32:245--256, 2012.

\bibitem{article2}
A.~Khan.
\newblock Bifurcation analysis of a discrete-time two-species model.
\newblock {\em Discrete Dynamics in Nature and Society}, 2020:1--12, 01 2020.

\bibitem{bib25}
Xiaojun Liu, Ling Hong, Lixin Yang, and Dafeng Tang.
\newblock Bifurcations of a new fractional-order system with a one-scroll
  chaotic attractor.
\newblock {\em Discrete Dynamics in Nature and Society}, 2019:8341514, Jan
  2019.

\bibitem{miles1983single}
Richard Miles and Robert~KS Wong.
\newblock Single neurones can initiate synchronized population discharge in the
  hippocampus.
\newblock {\em Nature}, 306(5941):371--373, 1983.

\bibitem{nagy2014efficient}
AM~Nagy and NH~Sweilam.
\newblock An efficient method for solving fractional hodgkin--huxley model.
\newblock {\em Physics Letters A}, 378(30-31):1980--1984, 2014.

\bibitem{bib30}
Paul Redfern and Stephen Thesleff.
\newblock Action potential generation in denervated rat skeletal muscle: Ii.
  the action of tetrodotoxin.
\newblock {\em Acta physiologica scandinavica}, 82(1):70--78, 1971.

\bibitem{article1}
M.~Sayed, Abd~allah Mousa, and Ibrahim Mustafa.
\newblock Stability and bifurcation analysis of a buckled beam via active
  control.
\newblock {\em Applied Mathematical Modelling}, 82, 02 2020.

\bibitem{bib1}
R.~D. Traub, R.~K. Wong, R.~Miles, and H.~Michelson.
\newblock A model of a ca3 hippocampal pyramidal neuron incorporating
  voltage-clamp data on intrinsic conductances.
\newblock {\em Journal of Neurophysiology}, 66(2):635--650, 1991.
\newblock PMID: 1663538.

\bibitem{bib27}
Seth~H Weinberg.
\newblock Membrane capacitive memory alters spiking in neurons described by the
  fractional-order hodgkin-huxley model.
\newblock {\em PloS one}, 10(5):e0126629, 2015.

\bibitem{bib26}
Svante Westerlund and Lars Ekstam.
\newblock Capacitor theory.
\newblock {\em IEEE Transactions on Dielectrics and Electrical Insulation},
  1(5):826--839, 1994.

\bibitem{PhysRevResearch.2.023281}
Can Xu, Xuebin Wang, and Per~Sebastian Skardal.
\newblock Bifurcation analysis and structural stability of simplicial
  oscillator populations.
\newblock {\em Phys. Rev. Research}, 2:023281, Jun 2020.

\bibitem{bib15}
Yue Zhang, Kuanquan Wang, Yongfeng Yuan, Dong Sui, and Henggui Zhang.
\newblock Stability and bifurcation analysis of hodgkin-huxley model.
\newblock In {\em 2013 IEEE International Conference on Bioinformatics and
  Biomedicine}, pages 49--54. IEEE, 2013.

\bibitem{bib12}
Jan Čermák and Luděk Nechvátal.
\newblock Local bifurcations and chaos in the fractional rössler system.
\newblock {\em International Journal of Bifurcation and Chaos}, 28(08):1850098,
  2018.

\end{thebibliography}

\end{document}